\documentclass[12pt]{article}

\usepackage{amsmath}
\usepackage{amssymb}
\usepackage[mathscr]{eucal}
\usepackage{graphicx}
\usepackage{color}

\usepackage{bm}
\usepackage{mathrsfs}
\usepackage{amsthm}
\usepackage{enumerate}
\usepackage[mathscr]{eucal}
\usepackage{eqlist}
\usepackage{graphicx}

\newtheorem{Theorem}{\sc Theorem}
\newtheorem{Definition}[Theorem]{\sc Definition}
\newtheorem{Proposition}[Theorem]{\sc Proposition}
\newtheorem{Lemma}[Theorem]{\sc Lemma}

\newtheorem{Remark}[Theorem]{\sc Remark}

\newtheorem{Problem}[Theorem]{\sc Problem}

\newcommand{\R}{{\if mm {\rm I}\mkern -3mu{\rm R}\else \leavevmode
		\hbox{I}\kern -.17em\hbox{R} \fi}}

\newcommand{\Div}{\mbox{\rm Div\,}}
\newcommand{\bu}{\mbox{\boldmath{$u$}}}
\newcommand{\bdelta}{\mbox{\boldmath{$\delta$}}}
\newcommand{\etab}{\mbox{\boldmath{$\eta$}}}
\newcommand{\bv}{\mbox{\boldmath{$v$}}}
\newcommand{\bw}{\mbox{\boldmath{$w$}}}
\newcommand{\bx}{\mbox{\boldmath{$x$}}}

\newcommand{\fb}{\mbox{\boldmath{$f$}}}

\newcommand{\bxi}{\mbox{\boldmath{$\xi$}}}

\newcommand{\bsigma}{\mbox{\boldmath{$\sigma$}}}
\newcommand{\btau}{\mbox{\boldmath{$\tau$}}}
\newcommand{\bvarepsilon}{\mbox{\boldmath{$\varepsilon$}}}
\newcommand{\bnu}{\mbox{\boldmath{$\nu$}}}

\newcommand{\bzero}{\mbox{\boldmath{$0$}}}

\def\sqr#1#2{{
		\vcenter{
			\vbox{\hrule height.#2pt
				\hbox{\vrule width.#2pt height#1pt \kern#1pt
					\vrule width.#2pt
				}
				\hrule height.#2pt
			}
		}
}}

\def\Div{\mathop{\rm Div}\nolimits}

\def\bar{\overline}
\def\real{\mathbb{R}}

\def\lista#1
{{ \itemindent 0.0cm \labelsep .2cm \leftmargin 0.8cm \rightmargin
		0.0cm \labelwidth 0.6cm \topsep 0.0mm
		\parsep 0.0mm
		\itemsep 0.0mm
		\begin{list}{}
			{ \setlength{\leftmargin}{.8cm} \setlength{\rightmargin}{0.0cm}
				\setlength{\parsep}{0.0mm} \setlength{\topsep}{.0mm}
				\setlength{\parskip}{.0cm} \setlength{\itemsep}{.0cm} }
			{#1}\end{list}} }

\begin{document}
\title{Rothe method and numerical analysis for history-dependent hemivariational inequalities with applications to contact mechanics 		
	\thanks{
	\, Project supported by the European Union's Horizon 2020 Research and Innovation Programme under the Marie Sk{\l}odowska-Curie grant agreement No. 823731 CONMECH, the National Science Center of Poland under Maestro Project No. UMO-2012/06/A/ST1/00262, and  National Science Center of Poland under Preludium Project No. 2017/25/N/ST1/00611. It is also supported by the International Project co-financed by the Ministry of Science and Higher Education of Republic of Poland under Grant No. 3792/GGPJ/H2020/2017/0.}}
	
\author{
	Stanis{\l}aw Mig\'orski\footnote{\,
		College of Applied Mathematics, Chengdu University of Information Technology, Chengdu,
        610225, Sichuan Province, P.R. China, and
		Jagiellonian University in Krakow, Chair of Optimization and Control, ul. Lojasiewicza 6, 30348 Krakow, Poland. E-mail address: stanislaw.migorski@uj.edu.pl. Tel.: +48-12-6646666.} \ \  and \ \
	Shengda Zeng\footnote{\,
		Jagiellonian University in Krakow, Faculty of Mathematics and Computer Science, ul. Lojasiewicza 6, 30348 Krakow, Poland. Corresponding author.
		E-mail address:	zengshengda@163.com; shengdazeng@gmail.com; shdzeng@hotmail.com.
Tel.: +86-18059034172.}
}
	
\date{}
\maketitle
	
\noindent {\bf Abstract.} \
In this paper an abstract evolutionary hemivariational inequality with a history-de\-pen\-dent operator is studied.
First, a result on its unique solvability and solution regularity is proved by applying the Rothe method.
Next, we introduce a numerical scheme to solve the inequality and derive error estimates.
We apply the results to a quasistatic frictional contact problem in which the material is modeled
with a viscoelastic constitutive law,
the contact is given in the form of multivalued normal compliance, and friction is described with a subgradient
of a locally Lipschitz potential. Finally, for the contact problem we provide the optimal error estimate.


\smallskip
	
	
\smallskip
	

\section{Introduction}\label{Introduction}

In this paper we are concerned with the existence and uniqueness of a solution to an abstract evolutionary hemivariational inequality which involve a history-dependent operator of the form
\begin{equation}\label{111}
\langle A u'(t)+Bu(t)+(\mathcal{R}u)(t) -f(t), v
\rangle+ J^0(Mu(t);Mv)\ge 0
\end{equation}

\noindent
for all $v \in V$, a.e. $t \in (0, T)$ with $u(0) = u_0$.
Here $A$ and $B$ are operators from a reflexive Banach space $V$ to its dual $V^*$, $M$ is a linear,
bounded operator, $J^0$ denotes the generalized gradient of a locally Lipschitz function, $f \colon (0, T) \to V^*$ and $u_0 \in V$ are given, and ${\cal R}$ represents a history-dependent operator.

The motivation to study the inequality of the form (\ref{111}) comes from contact problems in solid mechanics. It is known that when the external forces and tractions evolve slowly in time in such a way that the acceleration in the system is rather small and negligible, then the inertial
terms can be neglected. In such a way, we obtain
the quasistatic approximation (equilibrium equation) for the equation of motion.
Quasistatic contact models have been studied in several monographs and many papers dedicated to such phenomena,
see~\cite{DL,HS,SST,SHS} and the references therein.

In the first part of the paper, we deal with an abstract time-dependent hemivariational inequality of the form (\ref{111}).
The main results are delivered on existence, uniqueness and
regularity of a solution to the abstract hemivariational inequality, see Theorem~\ref{them3.2}.
We apply the Rothe method, see~\cite{Kac1,Kac2},
combined with a surjectivity result for a multivalued  and coercive operator.
The hemivariational inequality (\ref{111}) without a history-dependent operator has been recently investigated in~\cite{MOAN} by using the vanishing acceleration method, where a local existence result was proved.
In contrast to Theorem~17 of~\cite{MOAN}, here we provide a result on the global unique solvability to (\ref{111}). Also, our proof is now based on the Rothe method and is simpler, since we have eliminated the additional space $Z$ required in~\cite{MOAN}. Moreover, being motivated by applications to contact mechanics in Section~\ref{Dynamic}, the inequality (\ref{111}) involves a history-dependent operator.
We recall that the notion of a history-dependent operator
is quite recent and it was introduced in~\cite{SM2011}.
Various problems with history-dependent operators have been studied for the evolution variational and hemivariational inequalities
in~\cite{radulescu3,radulescu1,HMZ,HMS2017,MI2012,migorskizeng,migorskizeng2,migorskizeng3,MOSAMMA2,MOGORZ, MOGORZ2,radulescu2,zengliumigorski,zengmigorski}, and for the quasistatic problems
in~\cite{liumigorskizeng,MOS2015,SHM,SHMSUB,SP,SX,zengmigorski2}.
Furthermore, we study a fully discrete approximation for the problem (\ref{111}) which consists in finite difference discretization in time and finite element approximation in the spatial variable. We prove in Theorem~\ref{CEA} the C\'ea type error estimate for the hemivariational inequality.

In the second part of the paper, we apply the abstract results to a quasistatic frictional contact model for viscoelastic materials. The process is described by multivalued versions
of the nonmonotone normal compliance and friction boundary conditions.
We provide the variational formulation of the contact problem
for which we deliver a result on its unique weak global solvability.
In this way we improve the local existence result  of~\cite[Theorem~17]{MOAN}.
Finally, for the frictionless contact we establish a result on an optimal error estimate for the fully discrete approximation scheme.
Note that results on numerical anlaysis for hemivariational inequalities can be found in~\cite{Bartosz,HSB2017,HSS,HMP,SHM}
and the references therein.

The outline of the paper is as follows.
After recalling the basic notation in Section~\ref{Preliminaries}, in
Section~\ref{AbstractHVI} we formulate the abstract hemivariational inequality with a history-dependent operator.
In Section~\ref{Rothe} we apply the Rothe method to deliver existence and uniquence result for this inequality.
The error estimate of the C\'ea type for a fully discrete approximation is provided in Section~\ref{numerical}.
Finally, in Section~\ref{Dynamic}, we illustrate the
applicability of our results to the quasistatic frictional contact problem for viscoelastic material.

\section{Preliminaries}\label{Preliminaries}

In this section we recall the basic notation and some results which are needed in the sequel, see~\cite{cf,DMP1,DMP2,Zeidler}.
We use the standard notation for the Lebesgue and Sobolev spaces of functions defined on a finite time interval $[0,T]$ with values in a Banach space.
We denote by $\mathcal{L}(E,F)$
the space of linear and bounded operators
from a Banach space $E$ to a Banach space $F$ endowed with
the usual norm $\|\cdot\|_{\mathcal{L}(E, F)}$.
For a subset $S$ of Banach space $(E, \| \cdot \|_E)$,
we write
$\|S\|_E=\sup\{\|s \|_E \mid s \in S \}$.

%
Let $Y$ be a reflexive Banach space
and $\langle \cdot, \cdot \rangle$
denote the duality of $Y$ and $Y^*$.
A single-valued mapping $A\colon Y \to Y^*$ is called monotone if
$\langle Au-Av,u-v\rangle\ge 0$ for all
$u$, $v\in Y$.
An operator $A \colon Y \to Y^*$ is pseudomonotone
if for every sequence $\{ y_n \} \subseteq Y$ converging weakly to $y \in Y$ such that
$\displaystyle
\limsup \langle A y_n, y_n - y \rangle \le 0$,
we have
\begin{eqnarray*}
\displaystyle
\langle A y, y - z \rangle \le \liminf \langle A y_n, y_n - z \rangle \ \ \mbox{for all} \ \ z \in Y.
\end{eqnarray*}
Note that the operator $A \colon Y \to Y^*$
is pseudomonotone
if and only if the conditions
$y_n \to y$ weakly in $Y$ and
$\displaystyle
\limsup \langle A y_n, y_n - y \rangle \le 0$
entail
$\displaystyle
\lim \langle A y_n, y_n - y \rangle = 0$ and $A y_n \to A y$ weakly in $Y^*$.
It is also easy to check that if $A \in {\mathcal L}(Y, Y^*)$ is nonnegative, then it is pseudomonotone.

We recall the notion of the pseudomonotonicity
for a multivalued operator.

\begin{Definition}\label{DEFI3}
Let $Y$ be a reflexive Banach space.
An operator $T\colon Y\to 2^{Y^*}$ is pseudomonotone
if

\smallskip

\noindent
\rm{(a)} \ for every $v \in Y$, the set $Tv\subset Y^*$ is nonempty, closed and convex,

\smallskip

\noindent
\rm{(b)} \ $T$ is upper semicontinuous from each finite dimensional subspace of $Y$ to $Y^*$ endowed with the weak topology,

\smallskip

\noindent
\rm{(c)} \ for any sequences $\{u_n\}\subset Y$ and $\{u^*_n\}\subset Y^*$ such that $u_n\to u$ weakly in $Y$, $u_n^*\in Tu_n$ for all $n\geq 1$ and $\limsup\, \langle u_n^*,u_n-u\rangle\leq 0$, we have that for every $v\in Y$, there exists $u^*(v)\in Tu$ such that
    \begin{eqnarray*}
    \langle u^*(v),u-v\rangle\leq \liminf_{n\to\infty}\, \langle u_n^*,u_n-v\rangle.
    \end{eqnarray*}
\end{Definition}

We recall the following fundamental surjectivity theorem,
see~\cite[Theorem~1.3.70]{DMP2} or \cite{Zeidler}, which
will be used to prove existence of a solution to a static hemivariational inequality in Section~\ref{Rothe}.
\begin{Theorem}\label{surjective}
	Let $Y$ be a reflexive Banach space and
	$T \colon Y \to 2^{Y^*}$ be pseudomonotone and coercive. Then $T$ is surjective, i.e.,
	for every $f \in Y^*$, there is $u \in Y$
	such that $T u \ni f$.
\end{Theorem}

We hereafter recall the definition of the Clarke subgradient.
\begin{Definition}\label{SUB}
	Given a locally Lipschitz function
	$J \colon E \to \mathbb{R}$ on a Banach space $E$,
	we denote by $J^0 (u; v)$
	the generalized (Clarke) directional derivative
	of $J$ at the point $u\in E$ in the direction
	$v\in E$ defined by
	\begin{equation*}\label{defcalark}
	J^0(u;v) =
	\limsup
	\limits_{\lambda\to 0^{+}, \, w\to u}
	\frac{J(w+\lambda v)-J(w)}{\lambda}.
	\end{equation*}
	The generalized gradient of $J \colon E \to \mathbb{R}$ at $u\in E$ is defined by
	\begin{equation*}
	\partial J(u) =
	\{\, \xi\in E^{*} \mid J^0 (u; v)\ge
	\langle\xi, v\rangle
	\ \ \mbox{\rm for all} \ \ v \in E \, \}.
	\end{equation*}
\end{Definition}

The following result provides an example of a multivalued pseudomonotone operator which is a superposition of the Clarke subgradient with a compact operator.
The proof can be found in~\cite[Proposition 5.6]{Bartosz}.
\begin{Proposition}\label{proshvi2}
Let $V$ and $X$ be reflexive Banach spaces,  $M \colon V \to X$ be a linear, bounded, and compact operator.
We denote by $M^* \colon X^* \to V^*$
the adjoint operator of $M$.
Let $J\colon X \to \mathbb{R}$
be a locally Lipschitz function such that
\begin{equation*}
\| \partial J(v) \|_{X^{*}} \le c \, (1 + \| v \|_{X})
\ \ \mbox{\rm for all} \ \ v \in X
\end{equation*}
with $c>0$.
Then the multivalued operator
$F \colon V \to 2^{V^{*}}$
defined by
$F(v) = M^{*}\partial J(M v)$ for $v \in V$
is pseudomonotone.
\end{Proposition}

We conclude this section with a discrete version
of the Gronwall inequa\-lity
whose proof can be found in~\cite[Lemma 7.25]{HS}.
\begin{Lemma}\label{gronwall}
	Let $T>0$ be given. For a positive integer $N$, we define $\tau=\frac{T}{N}$.
	Assume that $\{g_{n}\}_{n=1}^{N}$ and $\{e_{n}\}_{n=1}^{N}$ are two sequences
	of nonnegative numbers satisfying
	\begin{equation*}
	e_{n} \le \overline{c}\, g_{n} + \overline{c} \, \tau
\, 	\sum_{j=1}^{n-1} e_{j}
	\ \ \mbox{\rm for} \ \ n=1,\ldots, N
	\end{equation*}
	for a positive constant $\overline{c}$ independent
	of $N$ (or $\tau$). Then there exists a positive constant $c$, independent of $N$ (or $\tau$) such that
	\begin{equation*}
	e_{n} \le c \, \Big(
	g_{n} + \tau
	\sum_{j=1}^{n-1} g_{j} \Big)
	\ \ \mbox{\rm for} \ \ n=1,\ldots, N.
	\end{equation*}
\end{Lemma}

\section{History-dependent hemivariational inequalities}
\label{AbstractHVI}

In this section we introduce a class of history-dependent hemivariational inequalities. This class will be studied in Section~\ref{Rothe} where the existence and uniqueness result for this class of inequalities will be provided.
A fully discrete approximation for the inequalities in this class will be discussed in Section~\ref{numerical}.

We use the following standard notation,
see~\cite{DMP1, DMP2,MOGORZ,Zeidler} for details.
Let $V \subset H \subset V^*$ be an evolution triple of spaces. Recall that this means that
$V$ is a reflexive and separable Banach space,
$H$ is a separable Hilbert space,
and the embedding $V \subset H$ is dense and continuous.
Let $i$ be the embedding operator between $V$ and $H$
which is assumed to be compact.
It is known that the adjoint operator
$i^* \colon H \to V^*$ is also linear, continuous and compact.
The duality pairing between $V^*$ and $V$, and a norm in $V$,
are denoted by $\langle\cdot,\cdot\rangle$ and $\|\cdot\|$,  respectively.
For the Hilbert space $H$, we denote its scalar product and a norm by $(\cdot,\cdot)$ and $\|\cdot\|_{H}$, respectively.

Given $0< T < +\infty$, let
$\mathcal{V} = L^2(0,T; V)$
and $\mathcal{H} = L^2(0,T; H)$.
It follows from the reflexivity of $V$ that both
$\mathcal{V}$ and its dual space
$\mathcal{V}^* = L^2(0,T; V^*)$
are reflexive Banach spaces as well.
Identifying $\mathcal{H} = L^2(0,T; H)$ with its dual,
we have the continuous embeddings
$\mathcal{V} \subset \mathcal{H} \subset \mathcal{V^*}$.

The notation
$\langle\cdot,\cdot
\rangle_{\mathcal{V}^{*}\times\mathcal{V}}$
stands for the duality pairing between $\mathcal{V}$ and $\mathcal{V}^{*}$.
Moreover, by $C(0, T; V)$ we denote the spave of continuous
functions on $[0, T]$ with values in $V$.


Let $X$ be a separable and reflexive Banach space.
Given operators
$A$, $B \colon V\to V^*$,
$M \colon V \to X$,
the function $J\colon X\to\real$,
$f\in \mathcal{V}^*$ and $u_{0}\in V$,
we consider the following evolutionary hemivariational inequality involving a history-dependent operator.
\begin{Problem}\label{p1}
	{\it Find an element $u \in {\mathcal V}$ such that $u' \in {\mathcal V}$  and }
	\begin{eqnarray}\label{hv}
		\left\{\begin{array}{lll}
			\big\langle A u'(t)+Bu(t)+(\mathcal{R}u)(t) -f(t), v
			\big\rangle+ J^0(Mu(t);Mv)\ge 0 \\ [2mm]
			\hspace{5cm}
			\mbox{\rm for all} \ \ v \in V, \, \mbox{\rm a.e.} \ \ t \in (0, T),
			\nonumber \\ [2mm]
			u(0)=u_{0}.\nonumber
		\end{array}\right.
	\end{eqnarray}
\end{Problem}
Here $\mathcal{R}\colon C(0,T;V)\to C(0,T;V^*)$ is an operator defined by
\begin{eqnarray}\label{defR}
(\mathcal{R}u)(t) =
E \, \bigg(\int_0^t q(t,s) u(s)\, ds + \alpha \bigg)
\ \ \mbox{for}\ \ t \in[0,T],
\end{eqnarray}
where $E\colon V\to V^*$, $\alpha\in V$ and
$q\colon [0,T]\times [0,T]\to \mathcal{L}(V,V)$.

\medskip

We impose the following assumptions on the data of Problem~\ref{p1}.

\medskip

\noindent
${\underline{H(A)}}$:\
The operator $A\colon V\to V^*$ is linear, bounded,
coercive and symmetric, i.e.,

\smallskip

\noindent

(i) $A\in\mathcal{L}(V,V^*)$.

\smallskip

\noindent

(ii) $\langle Av,v\rangle \ge m_A\|v\|^2$
for all $v\in V$ with $m_A>0$.

 \smallskip

\noindent

(iii) $\langle Av,w\rangle=\langle Aw,v\rangle$
for all $v$, $w\in V$.

\medskip

\noindent
${\underline{H(B)}}$:\ The operator $B\colon V\to V^*$ is linear, bounded and coercive, i.e.,

\smallskip

\noindent

(i) $B\in\mathcal{L}(V,V^*)$.

\smallskip

\noindent

(ii) $\langle Bv,v\rangle \ge m_B\|v\|^2$
for all $v\in V$ with $m_B>0$.

\medskip

\noindent
${\underline{H(E)}}$:\  $E\in\mathcal{L}(V,V^*)$.

\medskip

\noindent
${\underline{H(q)}}$:\  The function
$q\in C([0,T]\times[0,T], \mathcal{L}(V,V))$
is Lipschitz continuous with respect
to the first variable, i.e., there exists $L_q>0$
such that
\begin{eqnarray*}
\|q(t_1,s)-q(t_2,s)\|\le L_q|t_1-t_2|
\ \ \mbox{for all} \ \ t_1, t_2, s\in[0,T].
\end{eqnarray*}

\medskip

\noindent
${\underline{H(J)}}$:\ The functional $J\colon X\to \real$ is such that

\smallskip

\noindent

(i) $J$ is locally Lipschitz.

\smallskip

\noindent

(ii) There exists $c_J>0$ such that
$\|\partial J(u) \|_{X^*}\le c_J(1+\|u\|_X)$ for all $u\in X$.

\smallskip

\noindent

(iii) There exists $m_J \ge 0$ such that
\begin{eqnarray*}
\langle \xi-\eta,u-v\rangle_{X^*\times X}
\ge -m_J\|u-v\|^2_X,
\end{eqnarray*}
for all $u$, $v\in X$ and
$\xi\in\partial J(u)$, $\eta\in\partial J(v)$.

\medskip

\noindent
${\underline{H(f)}}$:\
$f\in \mathcal{V}^*$.

\medskip

\noindent
${\underline{H(M)}}$:\ The operator $M\colon V\to X$ is linear, continuous and compact.

\medskip

\noindent
${\underline{(H_0)}}$:\ $m_B>m_J\|M\|^2$.

\medskip

\begin{Remark}\label{remarkrelax}
Hypothesis $H(J)$(iii) is called the relaxed monotonicity condition for a locally Lipschitz function $J$. It was used in the literature (cf. \cite[Section 3.3]{smo1}) to ensure the uniqueness of the solution to hemivariational inequalities. This hypothesis has the equivalent formulation as follows
\begin{eqnarray*}
J^0(u;v-u)+J^0(v;u-v)\le m_J\|u-v\|_X^2,
\end{eqnarray*}
for all $u$, $v\in X$. In addition, examples of nonconvex functions which satisfy the relaxed monotonicity condition can be found in \cite{smo1,SHM}. Particularly, it can be proved that for a convex function, condition $H(J)$(iii) holds with $m_J=0$.
\end{Remark}

We recall, cf.~\cite{SM2011}, that an operator
${\cal S} \colon C(0, T; V) \to C(0, T; V^*)$
is called a history-dependent operator if there exists $L > 0$ such that
\begin{equation}\label{HIST}
\| ({\cal S} u_1)(t) - ({\cal S} u_2)(t) \|_{V^*}
\le L \int_0^t \| u_1 (s) - u_2(s) \|_V \, ds
\end{equation}
for all $u_1$, $u_2 \in C(0, T; V)$ and all $t \in [0, T]$.
We remark that under hypotheses $H(E)$, $H(q)$
and $\alpha \in V$,
the operator $\mathcal R$ defined in (\ref{defR})
satisfies condition (\ref{HIST})
with $L =c_Ec_q$, where $c_E =\|E\|$
and $c_q =\max_{(t,s)\in[0,T]\times[0,T]}\|q(t,s)\|$.

\section{Rothe method}\label{Rothe}

In this section, we present a result on existence and uniqueness of solution for Problem~\ref{p1}. The technique of proof relies on the Rothe method (known also as a method of lines, see~\cite{Kac1,Kac2}).
It consists in a time discretization in which we define an approximate sequence of functions by using the implicit (backward) Euler formula. Next, in each time step,
we will solve a stationary hemivational inequality. Finally, we construct the piecewise constant and piecewise affine
interpolants and prove a convergence result.

In the rest of the section, we denote by $C>0$ a constant whose value may change from line to line.

\medskip

Let $N\in\mathbb{N}$ be fixed and denote
$f_{\tau}^{k} =
\frac{1}{\tau}\int_{t_{k-1}}^{t_{k}}f(s)\, ds$
 for $k=1,\ldots,N$, where $t_{k}=k\tau$
and $\tau=\frac{T}{N}$.
Now, we
discuss the following discretized problem called
the Rothe problem.

\medskip

\begin{Problem}\label{ROTHE1}
Find $\{u_{\tau}^{k}\}_{k=0}^{N}\subset V$
such that $u_{\tau}^{0}=u_{0}$
and
\begin{eqnarray}\label{equation3}
\displaystyle\langle A u_\tau^k+\tau Bu_\tau^k+\tau x_\tau^k,v\rangle+\tau\langle \partial J(Mu_\tau^k),Mv\rangle_{X^*\times X}
\displaystyle
\ni\langle \tau f_\tau^k+Au_\tau^{k-1},v\rangle,
\end{eqnarray}
for all $v\in V$ and for $k=1,2,\ldots, N$,
where $x_\tau^k \in V^*$ is defined by
\begin{eqnarray*}
x_\tau^k =
E\bigg(\alpha+\sum_{j=1}^k\int_{t_{j-1}}^{t_j}q(t_k,s) u_\tau^j \, ds\bigg).
\end{eqnarray*}
\end{Problem}

First, we shall prove the existence and uniqueness of a solution to Problem~\ref{ROTHE1}.
\begin{Lemma}\label{lem3.1}
Assume that $u_0\in V$, $H(A)$, $H(B)$, $H(E)$, $H(q)$,  $H(J)$, $H(M)$ and $(H_0)$ hold.
Then there exists $\tau_{0}>0$ such that, for all $\tau\in(0,\tau_{0})$, Problem~\ref{ROTHE1}
has a unique solution.
\end{Lemma}

\noindent
{\bf Proof.}
Let $u_{\tau}^{0}$, $u_{\tau}^{1},\ldots,u_{\tau}^{k-1}$
be given. We will prove that there exists a unique element $u_{\tau}^{k}\in V$ which satisfies inclusion~(\ref{equation3}).
To end this, we apply Theorem~\ref{surjective} to show that the operator $L\colon V\to 2^{V^*}$ defined by
\begin{eqnarray*}
Lv = Av+\tau Bv+\tau E\bigg(\int_{t_{k-1}}^{t_k}q(t_k,s)v \, ds\bigg) + \tau M^*\partial J(Mv)
\end{eqnarray*}
for all $v\in V$, is surjective.

First, we show that
there exists $\tau_{0}>0$ such that, for all $\tau\in(0,\tau_{0})$,
$L$ is a pseudomonotone operator.
Indeed, by hypotheses
$H(A)$(i)-(ii),
$H(B)$(i)-(ii),
$H(E)$ and $H(q)$,
we can easily get that the operator
\begin{equation}\label{OPERAL}
v\mapsto Av+\tau Bv+\tau E\bigg(\int_{t_{k-1}}^{t_k}q(t_k,s)v\,ds \bigg)
\end{equation}
is bounded, continuous and monotone for
$\tau\in(0,\tau_0)$, where
$\tau_0=\frac{m_B}{c_Ec_q}$ with
$c_E=\|E\|$ and $c_q=\max_{(t,s)\in[0,T]\times[0,T]}\|q(t,s)\|$.
From~\cite[Theorem 3.69]{smo1}, we conclude that
the operator defined by (\ref{OPERAL})
is pseudomonotone. On the other hand, taking into account assumptions $H(J)$(i)-(ii), $H(M)$ and
Proposition~\ref{proshvi2}, it is clear that
the operator $v\mapsto M^*\partial J(Mv)$
is pseudomonotone as well.
Therefore, by using~\cite[Proposition 3.59(ii)]{smo1},
we infer that $L$ is a pseudomonotone operator too.

Subsequently, we prove that the operator $L$ is coercive.
From hypothesis $H(J)$ we derive the estimate (see~\cite{HSB2017})
\begin{eqnarray*}
\langle \partial J(u),u\rangle_{X^*\times X}\ge -m_J\|u\|_X^2-c_J\|u\|_X
\end{eqnarray*}
for all $u\in X$.
This inequality together with
$H(A)$(ii), $H(B)$(ii), $H(E)$ and $H(q)$ implies
\begin{eqnarray*}
&&\bigg\langle A u+\tau Bu+\tau E\bigg(\int_{t_{k-1}}^{t_k}q(t_k,s)u\, ds\bigg)
+\tau M^*\partial J(Mu),u\bigg\rangle\nonumber\\[2mm]
&&
\quad\ge m_A \|u\|^2+\tau m_B\|u\|^2-\tau^2 c_Ec_q\|u\|^2
-\tau m_J\|M\|^2\|u\|^2-\tau c_J\|M\|\|u\| \nonumber\\[2mm]
&&
\qquad\ge\big(m_A+\tau(m_B-m_J\|M\|^2-\tau c_Ec_q)\big)\|u\|^2-\tau c_J\|M\|\|u\|
\end{eqnarray*}
for all $u\in V$.
From the smallness condition $(H_0)$,
we choose $\tau_0= \frac{m_B-m_J\|M\|^2}{c_Ec_q} > 0$.
Hence, we deduce that the operator $L$ is coercive
for all $\tau\in(0,\tau_0)$.
Therefore, by the use of Theorem~\ref{surjective}, we obtain that $L$ is surjective, i.e.,
Problem~\ref{ROTHE1} has at least one solution
$u_\tau^k\in V$.

For uniqueness part, we assume that
$\overline{u}_\tau^k$ and $\widetilde{u}_\tau^k$ are two solutions in $V$ of Problem~\ref{ROTHE1}, that is,
\begin{eqnarray*}
\displaystyle\langle A \overline{u}_\tau^k+\tau B\overline{u}_\tau^k+\tau\overline{x}_\tau^k+\tau M^*\partial J(M\overline{u}_\tau^k),v\rangle
\displaystyle
\ge\langle \tau f_\tau^k+Au_\tau^{k-1},v\rangle
\ \ \mbox{\rm for all} \ \ v\in V
\end{eqnarray*}
and
\begin{eqnarray*}
\displaystyle\langle A \widetilde{u}_\tau^k+\tau B\widetilde{u}_\tau^k+\tau\widetilde{x}_\tau^k+\tau M^*\partial J(M\widetilde{u}_\tau^k),v\rangle
\displaystyle
\ge\langle \tau f_\tau^k+Au_\tau^{k-1},v\rangle
\ \ \mbox{\rm for all} \ \ v\in V,
\end{eqnarray*}
where the elements
$\overline{x}_\tau^k$ and $\widetilde{x}_\tau^k$
are defined by
\begin{eqnarray*}
\overline{x}_\tau^k
=E\bigg(\alpha+\sum_{j=1}^{k-1}
\int_{t_{j-1}}^{t_j}q(t_k,s)u_\tau^j\, ds
+\int_{t_{k-1}}^{t_k}q(t_k,s)\overline{u}_\tau^k\, ds\bigg)
\end{eqnarray*}
and
\begin{eqnarray*}
\widetilde{x}_\tau^k
=E\bigg(\alpha+\sum_{j=1}^{k-1}
\int_{t_{j-1}}^{t_j}q(t_k,s)u_\tau^j\, ds
+\int_{t_{k-1}}^{t_k}q(t_k,s)\widetilde{u}_\tau^k\, ds\bigg),
\end{eqnarray*}
respectively.
We take $v= \widetilde{u}_\tau^k-\overline{u}_\tau^k$
in the first inequality and
$v= \overline{u}_\tau^k-\widetilde{u}_\tau^k$
in the second one.
We add the resulting inequalities to get
\begin{eqnarray*}
&&\langle A \overline{u}_\tau^k-A\widetilde{u}_\tau^k,\overline{u}_\tau^k-\widetilde{u}_\tau^k\rangle
+\tau\langle B \overline{u}_\tau^k-B\widetilde{u}_\tau^k,\overline{u}_\tau^k-\widetilde{u}_\tau^k\rangle
+\tau\langle \overline{x}_\tau^k-\widetilde{x}_\tau^k,\overline{u}_\tau^k-\widetilde{u}_\tau^k\rangle\\[2mm]
&&\quad + \, \tau\langle \partial J(M\overline{u}_\tau^k)-\partial J(M\widetilde{u}_\tau^k),M\overline{u}_\tau^{k}-M\widetilde{u}_\tau^{k}\rangle_{X^*\times X}\le 0.
\end{eqnarray*}
Hence
\begin{eqnarray*}
\big(m_A+\tau(m_B-m_J\|M\|^2-\tau c_Ec_q)\big)\|\overline{u}_\tau^k-\widetilde{u}_\tau^k\|^2\le0.
\end{eqnarray*}
The smallness condition $(H_0)$ guarantees that  $\overline{u}_\tau^k=\widetilde{u}_\tau^k$,
which completes the proof of this lemma.
\hfill$\Box$

\medskip

Next, we establish the estimates
for the solution of Problem~\ref{ROTHE1}.

\begin{Lemma}\label{lem3.2}
Under assumptions of Lemma~\ref{lem3.1},
there exists $\tau_{0}>0$ and $C>0$ independent of $\tau$, such that for all $\tau\in(0,\tau_{0})$, the solution
$\{u_{\tau}^{k}\}_{k=0}^{N}\subset V$
of Problem~\ref{ROTHE1} satisfy
\begin{eqnarray}
&&\max\limits_{k=1,2,\ldots,N}\|u_{\tau}^{k}\|\leq C,\label{equation4}\\
&&\sum_{k=1}^N\|u_{\tau}^{k}-u_{\tau}^{k-1}\|^2\leq C,\label{equation8}\\
&&\max\limits_{k=1,2,\ldots,N}\|\xi_{\tau}^{k}\|_{X^*}\leq C,\label{equation9}\\
&&\tau\sum_{k=}^N\bigg\|\frac{u_{\tau}^{k}-u_\tau^{k-1}}{\tau}\bigg\|^2\leq C,\label{equation10}
\end{eqnarray}
where $\xi_\tau^k\in\partial J(Mu_\tau^k)$.
\end{Lemma}

\noindent
{\bf Proof.}
We choose $v=u_\tau^k$ in (\ref{equation3}), then use the hypotheses $H(A)$, $H(B)$ and the equality
$$
2 \langle Au_{\tau}^{k}-Au_{\tau}^{k-1},u_{\tau}^{k} \rangle=
\langle Au_{\tau}^{k},u_{\tau}^{k}\rangle-\langle Au_{\tau}^{k-1},u_{\tau}^{k-1}\rangle
+\langle A(u_{\tau}^{k}-u_{\tau}^{k-1}),u_{\tau}^{k}-u_{\tau}^{k-1}
\rangle
$$
to get
\begin{eqnarray}\label{equation5}
&&
\frac{1}{2}\langle Au_{\tau}^{k},u_{\tau}^{k}\rangle-\frac{1}{2}\langle Au_{\tau}^{k-1},u_{\tau}^{k-1}\rangle +\frac{1}{2}\langle A(u_{\tau}^{k}-u_{\tau}^{k-1}),u_{\tau}^{k}-u_{\tau}^{k-1}
\rangle+\tau m_B\|u_\tau^k\|^2 \nonumber \\[2mm]
&&\quad-\tau m_J\|M\|^2\|u_\tau^k\|^2-c_J\|M\|\|u_\tau^k\|\le \tau\|x_\tau^k\|_{V^*}\|u_\tau^k\|+\tau\| f_\tau^k\|_{V^*}\|u_\tau^k\|.
\end{eqnarray}
Next, the assumptions $H(E)$ and $H(q)$
imply
\begin{eqnarray}\label{equation6}
\tau\|x_\tau^k\|_{V^*}\|u_\tau^k\|\le \tau^2c_Ec_q\sum_{j=1}^k\|u_\tau^j\|\|u_\tau^k\|+\tau c_E\|\alpha\|\|u_\tau^k\|.
\end{eqnarray}
Combining (\ref{equation5}) and (\ref{equation6}),
and using the Cauchy inequality with $\varepsilon>0$,
we have
\begin{eqnarray*}
&&
\frac{1}{2}\langle Au_{\tau}^{k},u_{\tau}^{k}\rangle-\frac{1}{2}\langle Au_{\tau}^{k-1},u_{\tau}^{k-1}\rangle+\tau\big(m_B-m_J\|M\|^2-\varepsilon-\tau c_Ec_q\big)\|u_\tau^k\|^2 \nonumber\\[2mm]
 &&\quad+\frac{m_A}{2}\|u_\tau^k-u_\tau^{k-1}\|^2
\le C\tau\bigg(\tau\sum_{j=1}^{k-1}\|u_\tau^j\|^2+\|f_\tau^k\|_{V^*}+1\bigg).
\end{eqnarray*}
We now choose $\varepsilon =\frac{m_B-m_J\|M\|^2}{2}$ and $\tau_0=\frac{\varepsilon}{c_Ec_q}$.
Then, for all $\tau\in(0,\tau_0)$, it follows
\begin{eqnarray*}
&&
\frac{1}{2}\langle Au_{\tau}^{k},u_{\tau}^{k}\rangle-\frac{1}{2}\langle Au_{\tau}^{k-1},u_{\tau}^{k-1}\rangle+\frac{m_A}{2}\|u_\tau^k-u_\tau^{k-1}\|^2\\[2mm]
&&\qquad\qquad \le C\tau\bigg(\tau\sum_{j=1}^{k-1}\|u_\tau^j\|^2+\|f_\tau^k\|_{V^*}+1\bigg).
\end{eqnarray*}
Summing the above inequalities for $k=1,\ldots,n$,
where $1\le n\le N$, and then applying $H(A)$,
we get
\begin{eqnarray*}
&&
\frac{m_A}{2}\|u_\tau^k\|^2+\frac{m_A}{2}\sum_{k=1}^n\|u_\tau^k-u_\tau^{k-1}\|^2 \le C\bigg(\tau\sum_{k=1}^{n-1}\|u_\tau^k\|^2+1\bigg).
\end{eqnarray*}
Now, we use the discrete version
of the Gronwall inequality in Lemma~\ref{gronwall},
to verify estimates (\ref{equation4}) and (\ref{equation8}).
The estimate (\ref{equation9}) follows directly from
(\ref{equation4}) and $H(J)$(ii).

\medskip

Denote $v_\tau^k=\frac{u_\tau^k-u_\tau^{k-1}}{\tau}$ for $k=1,\ldots,N$.
We take $v=-v_\tau^k$ in (\ref{equation3}) to get
\begin{eqnarray*}
&&m_A\|v_\tau^k\|^2-\|B\|\|u_\tau^k\|\|v_\tau^k\|-\|\xi_\tau^k\|_{X^*}\|M\|\|v_\tau^k\| -\|x_\tau^k\|_{V^*}\|v_\tau^k\|
\\[2mm]
&&\quad\le\bigg\langle\frac{Au_\tau^k-Au_\tau^{k-1}}{\tau},v_\tau^k\bigg\rangle+\langle x_\tau^k,v_\tau^k\rangle+\langle Bu_\tau^k,v_\tau^k\rangle+\langle\xi_\tau^k, Mv_\tau^k\rangle_{X^*\times X}
\\[2mm]
&&\qquad\le\langle f_\tau^k, v_\tau^k\rangle\le \|f_\tau^k\|_{V^*}\|v_\tau^k\|,
\end{eqnarray*}
hence,
\begin{eqnarray*}
&&m_A\|v_\tau^k\|^2\le \big(\|B\|\|u_\tau^k\|+\|\xi_\tau^k\|_{X^*}\|M\|+\|x_\tau^k\|_{V^*}+\|f_\tau^k\|_{V^*}\big)\|v_\tau^k\|
\end{eqnarray*}
The latter together with (\ref{equation4}),
(\ref{equation9}), $H(E)$, $H(q)$ and
the Cauchy inequality with $\varepsilon>0$ implies
\begin{equation*}
(m_A-\varepsilon)\, \|v_\tau^k\|^2\le C(1+\|f_\tau^k\|_{V^*}).
\end{equation*}
We choose now $\varepsilon=\frac{m_A}{2}$ to get
\begin{equation*}
\tau\sum_{k=1}^N\|v_\tau^k\|^2\le C(1+\tau \sum_{k=1}^N\|f_\tau^k\|_{V^*})\le C.
\end{equation*}
So, we obtain the estimate (\ref{equation10}), which completes the proof of this lemma.
\hfill$\Box$

\medskip

Subsequently, for a given $\tau>0$,
we define the piecewise affine function $u_{\tau}$ and the piecewise constant interpolant functions $\overline{u}_{\tau}$, $\xi_{\tau}$, ${f}_{\tau}$  and $w_\tau$ as follows
\begin{eqnarray}
&&u_{\tau}(t)=u_{\tau}^{k}+\frac{t-t_{k}}{\tau}(u_{\tau}^{k}-u_{\tau}^{k-1})\ \ \textrm{ for }\ \ t\in(t_{k-1},t_{k}],\nonumber\\[2mm]
&&\overline{u}_{\tau}(t)=\left\{\begin{array}{lll}
u_{\tau}^{k}, & t\in(t_{k-1},t_{k}],\nonumber\\[2mm]
u_{\tau}^{0},& t=0,
\end{array}\right.\nonumber\\[2mm]
&&\xi_\tau(t)=\xi_{\tau}^{k}, \ \  t\in(t_{k-1},t_{k}],\nonumber\\[2mm]
&&{f}_{\tau}(t)=\left\{\begin{array}{lll}
f_{\tau}^{k}, & t\in(t_{k-1},t_{k}],\nonumber\\[2mm]
f(0), & t=0,
\end{array}\right.\nonumber\\[2mm]
&&w_\tau(t)=
\left\{\begin{array}{lll}
\alpha+\sum_{j=1}^k\int_{t_{j-1}}^{t_j}q(t_k,s)u_\tau^j\,ds,
&t\in(t_{k-1},t_k],\nonumber\\[2mm]
\alpha,&t=0.\nonumber
\end{array}\right.
\end{eqnarray}

\medskip

Now, we rewrite Problem~\ref{ROTHE1}
in the following equivalent form
\begin{eqnarray}\label{integralform}
&&\langle Au_\tau'(t)+B\overline{u}_\tau(t)+E(w_\tau(t)),v\rangle +\langle\xi_\tau(t),Mv\rangle_{X^*\times X}=\langle f_\tau(t),v\rangle
\end{eqnarray}
for all $v\in V$ and a.e. $t\in(0,T)$, where $\xi_\tau(t)\in \partial J(M\overline{u}_\tau(t))$  for a.e. $t\in(0,T)$.

\medskip

The main results of this section is delivered
in the following theorem.
\begin{Theorem}\label{them3.2}
Under assumptions of Lemma~\ref{lem3.1}, Problem~\ref{p1}
has a unique solution $u\in H^1(0,T;V)$.
\end{Theorem}

\noindent
{\bf Proof.}
The bound (\ref{equation4}) ensures that  $\{\overline{u}_\tau\}$ is bounded  in $\mathcal{V}$
due to the following inequality
\begin{eqnarray*}
\|\overline{u}_\tau\|^2_{\mathcal{V}}=\tau\sum_{n=1}^N\|u_\tau^n\|^2\le C.
\end{eqnarray*}
It follows from the reflexivity of $\mathcal{V}$
that there exists a function $u\in \mathcal{V}$ such that, passing to a subsequence again indexed by $\tau$, we have
\begin{eqnarray}\label{equation12}
\overline{u}_\tau\to u\ \, \mbox{ weakly in }\mathcal{V},\ \, \mbox{as }\tau\to 0.
\end{eqnarray}
Also, from (\ref{equation4}), we have that the sequence $\{u_\tau\}$ is bounded in $\mathcal{V}$, and therefore, there exists $u_1\in\mathcal{V}$ such that
\begin{eqnarray}\label{equation13}
u_\tau\to u_1\ \, \mbox{ weakly in }\mathcal{V},\ \, \mbox{as }\tau\to 0.
\end{eqnarray}
Hence, we get $\overline{u}_\tau-u_\tau \to u-u_1$ weakly in $\mathcal{V}$, as $\tau\to 0$.
By the H\"older inequality and the boundedness of $\{u_\tau'\}$ (see (\ref{equation10}))
\begin{eqnarray*}
\|u_\tau'\|^2_{\mathcal{V}}=\tau\sum_{k=1}^N\|v_\tau^k\|^2\le C,
\end{eqnarray*}
we have
\begin{eqnarray}\label{equation14}
&&
\|\overline{u}_\tau-u_\tau\|_{\mathcal{V}}^2
=\sum_{k=1}^N\int_{t_{k-1}}^{t_k}(t_k-s)^2
\|v_\tau^k\|^2 \, ds \nonumber \\
&&\quad= \sum_{k=1}^N\int_{t_{k-1}}^{t_k}(t_k-s)^2
\|u_\tau'(s)\|^2 \, ds
\le\frac{\tau^2}{3} \|u_\tau'\|^2_{\mathcal{V}}.
\end{eqnarray}
From estimate (\ref{equation14}) we deduce that
$u=u_1$.
On the other hand, by the boundedness of $\{u_\tau'\}$ (see (\ref{equation10})), we also obtain (cf. \cite[Proposition~23.19. p. 419]{Zeidler})
\begin{eqnarray}\label{equation15}
u_\tau'\to u' \ \, \mbox{ weakly in }\mathcal{V},\ \, \mbox{as }\tau\to 0.
\end{eqnarray}
In addition, using the boundedness
of $\{\xi_\tau\}$ (see (\ref{equation9})) and the reflexivity of the space $\mathcal{X}^*$,
we conclude
\begin{eqnarray}\label{equation16}
\xi_\tau\to \xi \ \, \mbox{ weakly in }\mathcal{X}^*,
\ \, \mbox{as }\tau\to0 \ \ \mbox{with} \ \ \xi\in\mathcal{X}^*.
\end{eqnarray}
By virtue of the hypothesis $H(q)$ and boundedness of $\{\overline{u}_\tau\}$ (see (\ref{equation4})),
one has the following estimate
for $t\in(t_{k-1},t_k]$
\begin{eqnarray}\label{equation17}
&&
\bigg\|\int_0^tq(t,s)\overline{u}_\tau(s)\, ds
-\int_0^{t_k}q(t_k,s)\overline{u}_\tau(s)\, ds\bigg\|
\le \int_t^{t_k}\|q(t_k,s)\overline{u}_\tau(s)\|\, ds
\nonumber\\
&&\quad
+\int_0^t
\big\|(q(t,s)-q(t_k,s))\overline{u}_\tau(s)\big\|\,ds
\le C_0\tau
\end{eqnarray}
for some $C_0 > 0$, which is independent of $\tau$.
Moreover, \cite[Lemma 3.3]{Carstensen} implies that
\begin{eqnarray}\label{equation18}
f_\tau\to f \ \, \mbox{strongly in } \mathcal{V}^*,
\ \, \mbox{ as }\tau \to 0.
\end{eqnarray}

\smallskip

Next, we shall show that $u$ is a solution of
Problem~\ref{p1}. To this end, we define
the Nemytskii operators $\mathcal{A}$,  $\mathcal{B}\colon \mathcal{V}\to\mathcal{V}^*$
by $(\mathcal{A}v)(t)=A(v(t))$ and $(\mathcal{B}v)(t)=B(v(t))$ for all $v\in\mathcal{V}$ and a.e. $t\in (0,T)$.
From hypotheses $H(A)$ and $H(B)$,
it is clear that $\mathcal{A}$ and $\mathcal{B}$ are both linear and bounded,
so they are also weakly continuous.
Thus from (\ref{equation15}) and (\ref{equation12})
we obtain
$\mathcal{A}u_\tau'\to \mathcal{A}u'$ and $\mathcal{B}\,\overline{u}_\tau\to \mathcal{B}u$
both weakly in $\mathcal{V}^*$, as $\tau\to0$,
i.e.,
\begin{equation}
\lim_{\tau\to0}\langle \mathcal{A}u_\tau',v\rangle_{\mathcal{V}^*\times \mathcal{V}}= \langle\mathcal{A}u', v\rangle_{\mathcal{V}^*\times \mathcal{V}} \ \ \mbox{and}\ \
\lim_{\tau\to0}\langle \mathcal{B}\, \overline{u}_\tau, v\rangle_{\mathcal{V}^*\times \mathcal{V}} =\langle \mathcal{B}u,v\rangle_{\mathcal{V}^*\times \mathcal{V}} \label{equation19}
\end{equation}
for all $v\in\mathcal{V}$.
Now, we consider the Nemitskii operators
$\mathcal{E}$,
$\mathcal{E}_2\colon \mathcal{V}\to \mathcal{V}^*$
by
\begin{equation*}
(\mathcal{E}v)(t) =
E\bigg(\int_0^tq(t,s)v(s)\, ds\bigg)
\ \ \mbox{and} \ \ (\mathcal{E}_2v)(t) = Ev(t)
\end{equation*}
for all $v\in\mathcal{V}$ and a.e. $t\in (0,T)$.
It is obvious that $\mathcal{E}$
is weakly continuous being bounded and linear.
From the convergence (\ref{equation12}), one has
\begin{equation*}
\lim_{\tau\to0}\langle\mathcal{E}\overline{u}_\tau,v\rangle_{\mathcal{V}^*\times \mathcal{V}}=\langle\mathcal{E}u,v\rangle_{\mathcal{V}^*\times \mathcal{V}}
\end{equation*}
for all $v\in\mathcal{V}$.
Next, from $H(E)$, $H(q)$ and (\ref{equation17}),
we have
\begin{equation*}
\mathcal{E}_2(w_\tau-\alpha)
-\mathcal{E}(\overline{u}_\tau) \to 0
\ \, \mbox{strongly in}\ \mathcal{V}^*,
\ \, \mbox{ as }\tau \to 0
\end{equation*}
which implies
\begin{eqnarray}\label{equation20}
&&\lim_{\tau\to 0}\langle \mathcal{E}_2(w_\tau),
v\rangle_{\mathcal{V}^*\times{V}}\nonumber\\ &&\quad=\lim_{\tau\to 0}\big(\langle \mathcal{E}_2(w_\tau-\alpha)
-\mathcal{E}(\overline{u}_\tau),
v \rangle_{\mathcal{V}^*\times{\cal V}}
+\langle \mathcal{E}(\overline{u}_\tau),
v\rangle_{\mathcal{V}^*\times
	{\cal V}}+\langle \mathcal{E}_2(\alpha),
v\rangle_{\mathcal{V}^*\times
	{\cal V}}\big)\nonumber\\[2mm]
&&\qquad=\langle\mathcal{E}u,
v\rangle_{\mathcal{V}^*\times \mathcal{V}}+\langle \mathcal{E}_2(\alpha),
v\rangle_{\mathcal{V}^*\times
	{\cal V}}
\end{eqnarray}
for all $v\in\mathcal{V}$.

Since the embedding
$H^1(0, T; V) \subset C(0,T;V)$
is continuous, from the convergences (\ref{equation13}) and (\ref{equation15}),
by~\cite[Lemma 4(a)]{MOAN}, we have
\begin{equation}\label{equation22}
u_\tau(t)\to u(t)\ \,\mbox{ weakly in } V, \ \, \mbox{ as }\tau\to 0,
\end{equation}
for all $t\in[0,T]$.
Using the convergence
$\overline{u}_\tau - u_\tau\to 0$ strongly in $\mathcal{V}$, as $\tau\to 0$, by the converse Lebesgue dominated convergence theorem, \cite[Theorem 2.39]{smo1},
we may assume that $\overline{u}_\tau(t) - u_\tau(t)
\to 0$ strongly in $V$ for a.e. $t \in (0, T)$,
as $\tau\to 0$. This together with (\ref{equation22}) implies
\begin{equation*}
\overline{u}_\tau(t)\to u(t)\ \,\mbox{ weakly in }V, \ \, \mbox{for a.e. } t\in(0,T).
\end{equation*}
From the compactness of the operator
$M$, we deduce $M\overline{u}_\tau(t)\to Mu(t)$
strongly in $X$ for a.e. $t\in(0,T)$.
Since
$\xi_\tau(t)\in\partial J(M\overline{u}_\tau(t))$
for a.e. $t\in(0,T)$,
we use also the convergence (\ref{equation16}),
and by
\cite[Theorem 1, Section 1.4]{AC}, we have
\begin{equation}\label{equation23}
\xi(t)\in\partial J(Mu(t))
\ \,\mbox{ for a.e. }t\in(0,T).
\end{equation}
Now, we introduce the Nemitskii operator $\mathcal{M}\colon\mathcal{V} \to\mathcal{X}$
defined by $(\mathcal{M}v)(t)=M(v(t))$ for all $v\in\mathcal{V}$ and a.e. $t\in (0,T)$, so, from (\ref{equation16}), we have
\begin{equation}\label{equation24}
\lim_{\tau\to0}\langle \xi_\tau,\mathcal{M}v\rangle_{\mathcal{X}^*\times \mathcal{X}} =\langle \xi,\mathcal{M}v\rangle_{\mathcal{X}^*\times \mathcal{X}}
\end{equation}
for all $v\in\mathcal{V}$.

From (\ref{equation18})--(\ref{equation20}),  (\ref{equation23}) and (\ref{equation24}),
we infer that
\begin{equation}\label{gequality}
\langle \mathcal{A}u'+\mathcal{B}u+\mathcal{E}u+\mathcal{E}_2\alpha,v\rangle_{\mathcal{V}^*\times \mathcal{V}}+ \langle \xi,\mathcal{M}v\rangle_{\mathcal{X}^*\times \mathcal{X}}= \langle f,v\rangle_{\mathcal{V}^*\times \mathcal{V}}
\end{equation}
for all $v\in\mathcal{V}$ with $\xi(t)\in\partial J(Mu(t))$  for a.e. $t\in(0,T)$. Furthermore, we shall show that $u \in {\cal V}$ with $u' \in {\cal V}$ is also a solution of Problem~\ref{p1}.
Arguing by contradiction, we suppose that $u$ is not a solution to Problem~\ref{p1}. This means there exist a measurable set $I\subset [0,T]$ with $\mbox{meas(I)}>0$ and $v^*\in V$ such that
\begin{equation}\label{medinq}
\langle A u'(t)+Bu(t)+(\mathcal{R}u)(t), v^*
\rangle+ J^0(Mu(t);Mv^*)<\langle f(t), v^*
\rangle\ \,\mbox{ for a.e. }t\in I.
\end{equation}
We now denote a function $\widetilde{v}\in \mathcal V$ by
\begin{eqnarray*}
\widetilde v (t)=
\left\{\begin{array}{ll}
v^*&\mbox{ if }t\in I\\[2mm]
0&\mbox{otherwise}.
\end{array}\right.
\end{eqnarray*}
Inserting $v=\widetilde v$ into (\ref{gequality}) and taking account of (\ref{medinq}), it follows from~\cite[Theorem 3.47]{smo1} that
\begin{eqnarray*}
&&\int_I\langle f(t),v^*\rangle\,dt
\le \int_I\langle A u'(t)+Bu(t)+(\mathcal{R}u)(t) -f(t), v^*
\rangle+ J^0(Mu(t);Mv^*)\,dt
\\[2mm]
&&\qquad<\int_I\langle f(t),v^*\rangle\,dt.
\end{eqnarray*}
This results a contradiction, so, $u \in {\cal V}$ with $u' \in {\cal V}$ is also a solution of Problem~\ref{p1}.

\medskip

Finally, we will verify that the solution of
Problem~\ref{p1} is unique.
Let $u_1$ and $u_2$ be two solutions of
Problem~\ref{p1}. Then
\begin{equation*}
\big\langle A u'_1(t)+Bu_1(t)+(\mathcal{R}u_1)(t) -f(t), v
\big\rangle +J^0(Mu_1(t);Mv)\ge 0
\end{equation*}
and
\begin{equation*}
\big\langle A u'_2(t)+Bu_2(t)+(\mathcal{R}u_2)(t) -f(t), v
\big\rangle +J^0(Mu_2(t);Mv)\ge 0
\end{equation*}
for all $v \in V$ and a.e.  $t \in (0, T)$. Taking $v=u_2(t)-u_1(t)$ in the first inequality and $v=u_1(t)-u_2(t)$ in the second one, we add the resulting inequalities to get
\begin{eqnarray*}
&&\big\langle A u'_1(t)-A u'_2(t),  u_1(t)-u_2(t)
\big\rangle+\langle Bu_1(t)-Bu_2(t),u_1(t)-u_2(t)\rangle\nonumber\\[2mm]
&&\quad
\le J^0(Mu_1(t);Mv)+J^0(Mu_1(t);Mu) +\langle(\mathcal{R}u_1)(t)-(\mathcal{R}u_2)(t),u_2(t)-u_1(t)\rangle
\end{eqnarray*}
for a.e. $t\in(0,T)$. We use the assumptions $H(A)$, $H(B)$, $H(E)$, $H(q)$, and $H(J)$(iii) to obtain
\begin{eqnarray*}
&&\frac{1}{2}\frac{d}{dt}
\big\langle A (u_1(t)-u_2(t)), u_1(t)-u_2(t)
\big\rangle
+ (m_B-m_J\|M\|^2) \|u_1(t)-u_2(t)\|^2
\nonumber\\[2mm]
&&\quad\le c_Ec_q\int_0^t\|u_1(s)-u_2(s)\| \|u_1(t)-u_2(t)\| \, ds.
\end{eqnarray*}
We integrate this inequality on $[0,t]$, where $t\in[0,T]$, and use $H(A)$(ii) and $(H_0)$ to deduce
\begin{eqnarray*}
&&\frac{m_A}{2}\|u_1(t)-u_2(t)\|^2\le
\frac{1}{2}\big\langle A (u_1(t)-u_2(t)),
u_1(t)-u_2(t) \big \rangle\\
&&\quad
\le c_E\, c_q
\int_0^t\|u_1(s)-u_2(s)\|\int_0^s
\|u_1(\eta)-u_2(\eta)\|\, d\eta \, ds \\
&&\qquad\le c_E\, c_q\bigg(\int_0^t\|u_1(s)-u_2(s)
\| \, ds\bigg)^2
\end{eqnarray*}
for all $t\in[0,T]$. Hence
\begin{equation*}
\|u_1(t)-u_2(t)\|\le
\bigg(\frac{2c_E\,c_q}{m_A}\bigg)^{\frac{1}{2}}
\int_0^t\|u_1(s)-u_2(s)\| \, ds
\end{equation*}
for all $t\in[0,T]$. Finally, we use the Gronwall inequality (see e.g.~\cite[Lemma 2.31]{SHS})
to obtain $u_1=u_2$. This completes the proof of the theorem.
\hfill$\Box$

\section{A fully discrete approximation scheme}
\label{numerical}

In this section, we
study a fully discrete approximation scheme
for the history-dependent hemivariational inequality
stated in Problem~\ref{p1}.
In this method the time variable is discretized
by finite difference and the spatial variable is approximated by finite elements.

Assume that $V^h$ is a finite dimensional subspace of $V$ and $u_0^h\in V^h$ is an approximation of the initial point $u_0 \in V$.
For $N\in\mathbb{N}$, $N > 0$ given,
we denote the time step length by
$k=\frac{T}{N}$ and
$t_n=kn$ for $n=0,\ldots,N$.
For a continuous function $g$ defined on the interval $[0,T]$, in the sequel,  we will write
$g_n=g(t_n)$ for $n=0,\ldots,N$.
In addition, for a sequence $\{u_n\}_{n=0}^N$,
we use the notation
\begin{equation*}
\delta u_n=\frac{u_n-u_{n-1}}{k},\ \ n=1,\ldots,N.
\end{equation*}
For the history-dependent operator
\begin{equation*}
(\mathcal{R}v)(t)
= \int_0^tq(t,s)v(s)\, ds \ \,\mbox{for}\ \,
v\in C(0,T;V), \ t \in [0, T],
\end{equation*}
we introduce a modified trapezoidal approximation for  $\mathcal{R}$ defined by
\begin{equation}\label{trapezoidal}
\mathcal{R}_n^kv =
E\bigg(\sum_{j=1}^n\int_{t_{j-1}}^{t_j}q(t_n,s)v_j \,ds
+\alpha\bigg)
\end{equation}
for $v=\{v_j\}_{j=1}^N$.
In addition, if $w\in C(0,T;V)$, then the expression  $\mathcal{R}_n^kw$ is understood as follows
\begin{equation*}
\mathcal{R}_n^kw =
E\bigg(\sum_{j=1}^n\int_{t_{j-1}}^{t_j}q(t_n,s)w(t_j) \,ds
+\alpha\bigg).
\end{equation*}

Subsequently,
we consider the following fully discrete approximation problem for Problem~\ref{p1}.
\begin{Problem}\label{fullydiscrete}
Find $u^{hk} = \{u_n^{hk}\}\subset V^h$ such that $u_0^{hk}=u_0^h$ and
\begin{eqnarray}\label{equation26}
&&\langle A\delta u_n^{hk}+B u_n^{hk}+\mathcal{R}_n^{k}u^{hk},v^h-u_n^{hk}\rangle +J^0(Mu_n^{hk};Mv^{h}-Mu_n^{hk})\nonumber\\[2mm]
&&\hspace{3cm}\ge\langle f_n,v^h-u_n^{hk}\rangle
\ \ \mbox{for all} \ \ v^h\in V^h
\end{eqnarray}
for all $n=1,2,\ldots,N$.
\end{Problem}

We will provide an error analysis of the fully discrete approximation (\ref{equation26}).
Our goal is to prove the C\'ea type inequality
for Problem~\ref{fullydiscrete}.

First, exploiting the definition of $\delta u_n^{hk}$,
the inequality (\ref{equation26})
can be reformulated as follows
\begin{eqnarray}\label{equation27}
&&\langle A u_n^{hk}+kB u_n^{hk}+k\mathcal{R}_n^{k}u^{hk},v^h-u_n^{hk}\rangle +kJ^0(Mu_n^{hk};Mv^{h}-Mu_n^{hk})\nonumber\\[2mm]
&&\hspace{3cm}\ge\langle kf_n+Au_{n-1}^{hk},v^h-u_n^{hk}\rangle\
\ \ \mbox{for all} \ \  v^h\in V^h.
\end{eqnarray}
This inequality represents a stationary hemivariational inequality. When $k$ small enough, from Lemma~\ref{lem3.1}, we know
that under the hypotheses $H(A)$, $H(B)$, $H(E)$, $H(q)$, $H(J)$, $H(M)$ and $(H_0)$,
it has a unique solution $u_n^{hk}\in V^h$.
Moreover, Theorem~\ref{them3.2} reveals that
Problem~\ref{p1} has a unique solution
$u\in H^1(0,T;V)$.

Since $A\in \mathcal{L}(V,V^*)$ is coercive, in what follows, for a convenience, we introduce
the norm $\|\cdot\|_A$ by
$\|v\|_A^2=\langle Av,v\rangle$ for all $v \in V$,
which is equivalent to the norm $\|\cdot\|_V$.
In the sequel, we denote by $C>0$ a constant which may differ from line to line, but it is independent of $h$ and $k$.

For an error analysis, we have from (\ref{111})
at $t=t_n$ that
\begin{equation}\label{equation28}
\langle Au_n'+Bu_n+\mathcal{R}_nu, v-u_n\rangle+J^0(Mu_n;Mv-Mu_n)
\ge\langle f_n,v-u_n\rangle
\end{equation}
for all $v \in V$, where $\mathcal{R}_nu=(\mathcal{R}u)(t_n)$.
Denote the errors
\begin{equation*}
\delta_n =\delta u_n-u'_n
\ \ \, \mbox{ and }\ \ \, e_n =u_n-u_n^{hk}
\end{equation*}
for $n=1$, $2,\ldots,N$.
Taking $v=u_n^{hk}$ in (\ref{equation28}), one has
\begin{eqnarray}\label{equation29}
&&\langle A\delta u_n+Bu_n+\mathcal{R}_nu,u_n^{hk}-u_n\rangle+J^0(Mu_n;Mu_n^{hk}-Mu_n)\nonumber\\[2mm]
&&\hspace{3cm}
\ge\langle f_n,u_n^{hk}-u_n\rangle+\langle A\delta_n,u_n^{hk}-u_n\rangle.
\end{eqnarray}
We add (\ref{equation29}) and (\ref{equation26}) to get
\begin{eqnarray*}
&&\langle A\delta u_n+Bu_n+\mathcal{R}_nu,u_n^{hk}-u_n\rangle+\langle A\delta u_n^{hk}+B u_n^{hk}+\mathcal{R}_n^{k}u^{hk},v^h-u_n^{hk}\rangle\nonumber\\[2mm]
&&\quad+J^0(Mu_n;Mu_n^{hk}-Mu_n)+J^0(Mu_n^{hk};Mv^{h}-Mu_n^{hk})\nonumber\\[2mm]
&&\qquad\ge \langle f_n,v^h-u_n\rangle+\langle A\delta_n,u_n^{hk}-u_n\rangle
\end{eqnarray*}
for all $v^h\in V^h$.
Hence
\begin{eqnarray*}
&&\langle A\delta (u_n-u_n^{hk})+B(u_n-u_n^{hk}),u_n^{hk}-u_n\rangle
+\langle\mathcal{R}_nu-\mathcal{R}_n^{k}u^{hk},u_n^{hk}-u_n\rangle\nonumber\\[2mm]
&&\quad+\langle A\delta u_n^{hk}+B u_n^{hk}+\mathcal{R}_n^{k}u^{hk},v^h-u_n\rangle +J^0(Mu_n;Mu_n^{hk}-Mu_n)\nonumber\\[2mm]
&&\qquad+J^0(Mu_n^{hk};Mv^{h}-Mu_n^{hk})\ge \langle f_n,v^h-u_n\rangle+\langle A\delta_n,u_n^{hk}-u_n\rangle
\end{eqnarray*}
for all $v^h\in V^h$.
We use the fact that the function
$v\mapsto J^0(Mu;Mv)$ is subadditive (see e.g., \cite[Proposition 3.23(i)]{smo1}), to obtain
\begin{equation*}
J^0(Mu_n^{hk};Mv^{h}-Mu_n^{hk})\le J^0(Mu_n^{hk};Mv^{h}-Mu_n)+J^0(Mu_n^{hk};Mu_n-Mu_n^{hk}).
\end{equation*}
So, we have
\begin{eqnarray*}
&&\langle A\delta (u_n-u_n^{hk})+B(u_n-u_n^{hk}),u_n^{hk}-u_n\rangle
+\langle\mathcal{R}_nu-\mathcal{R}_n^{k}u^{hk},u_n^{hk}-u_n\rangle\nonumber\\[2mm]
&&\quad+\langle A\delta u_n^{hk}+B u_n^{hk}+\mathcal{R}_n^{k}u^{hk},v^h-u_n\rangle +J^0(Mu_n;Mu_n^{hk}-Mu_n)\nonumber\\[2mm]
&&\qquad+J^0(Mu_n^{hk};Mv^{h}-Mu_n)+J^0(Mu_n^{hk};Mu_n-Mu_n^{hk})\nonumber\\[2mm]
&&\qquad\quad\ge \langle f_n,v^h-u_n\rangle+\langle A\delta_n,u_n^{hk}-u_n\rangle
\end{eqnarray*}
for all $v^h\in V^h$. Combining this inequality with the identity
\begin{eqnarray*}
\langle A(u-v),u\rangle=\frac{1}{2}\big(\langle Au,u\rangle-\langle Av,v\rangle+\langle A(u-v),u-v\rangle\big)
\end{eqnarray*}
and using the hypotheses $H(B)$ and $H(J)$(iii)
(see Remark~\ref{remarkrelax}),
it follows that
\begin{eqnarray}\label{equation30}
&&\frac{1}{2k}\big(\|e_n\|^2_A-\|e_{n-1}\|^2_A\big)
+m_B\|e_n\|^2-m_J\|M\|^2\|e_n\|^2 -\|\mathcal{R}_nu-\mathcal{R}_n^{k}u^{hk}\|_{V^*}
\|e_n\|\nonumber\\[2mm]
&&\quad
\le\langle A\delta u_n^{hk}+B u_n^{hk}+\mathcal{R}_n^{k}u^{hk}-f_n,v^h-u_n\rangle
+J^0(Mu_n^{hk};Mv^{h}-Mu_n) \nonumber\\[2mm]
&&\qquad
-\langle A\delta_n,u_n^{hk}-u_n\rangle
\end{eqnarray}
for all $v^h\in V^h$. Furthermore,
we introduce a residual type quantity by
\begin{equation*}
S_n(v) =
\langle Au_n'+Bu_n+\mathcal{R}_nu-f_n,v-u_n\rangle
+J^0(Mu_n;Mv-Mu_n) \ \ \mbox{for} \ \ v\in V.
\end{equation*}
Using the fact that
$u\in H^1(0,T;V)$ (see Theorem~\ref{them3.2}),
we have
\begin{eqnarray*}
&&\|\mathcal{R}_nu-\mathcal{R}_n^ku\|_{V^*}\le c_E\sum_{i=1}^n\int_{t_{i-1}}^{t_i}
\big\|q(t_n,s)(u(s)-u_i)\big\|\, ds\nonumber\\
&&\quad\le C_1 k
\sum_{i=1}^n\int_{t_{i-1}}^{t_i}\|u'(s)\|\, ds=
C_1 k
\int_{0}^T\|u'(s)\|\, ds \le C_1\sqrt{T}\|u'\|_{\mathcal{V}}\,k \nonumber
\end{eqnarray*}
with some $C_1>0$, and
\begin{eqnarray*}
\|\mathcal{R}_n^ku-\mathcal{R}_n^{k}u^{hk}\|_{V^*}\le c_E\, c_q
\sum_{i=1}^n\int_{t_{i-1}}^{t_i}\|u_i-u_i^{hk}\|\, ds
\le k c_E c_q\sum_{i=1}^n\|u_i-u_i^{hk}\|.
\end{eqnarray*}
From these inequalities,
we obtain
\begin{equation}\label{equation31}
\|\mathcal{R}_nu-\mathcal{R}_n^{k}u^{hk}\|_{V^*}\le \|\mathcal{R}_nu-\mathcal{R}_n^ku\|_{V^*}
+\|\mathcal{R}_n^ku-\mathcal{R}_n^{k}u^{hk}\|_{V^*}
\le C_2 k \bigg( 1+\sum_{j=1}^n\|e_j\|\bigg),
\end{equation}
where $C_2=\max\{c_E c_q,C_1\sqrt{T}\|u'\|_{\mathcal{V}}\}$.
Therefore, from (\ref{equation30}), we have
\begin{eqnarray}\label{equation32}
&&\hspace{-1.0cm}
\frac{1}{2k}\big(\|e_n\|^2_A-\|e_{n-1}\|^2_A\big)- C_2k\bigg( 1+\sum_{j=1}^n\|e_j\|\bigg)\|e_n\|+(m_B-m_J\|M\|^2)
\|e_n\|^2
\nonumber\\
&&
\hspace{-1.0cm}
\quad\le \frac{1}{k}\langle Ae_n-Ae_{n-1},u_n-v_n^h\rangle+\langle Be_n,u_n-v_n^h\rangle+\langle \mathcal{R}_nu-\mathcal{R}_n^ku^{hk},u_n-v^h_n\rangle
\nonumber\\[2mm]
&&\hspace{-1.0cm}
\qquad+\langle\xi_n-\xi_n^{hk},
M(u_n-v^h_n)\rangle_{X^*\times X}
+S_n(v^h_n)+\langle A\delta_n,v^h_n-u_n\rangle
+\langle A\delta_n, e_n\rangle,
\end{eqnarray}
where $\xi_n\in\partial J(Mu_n)$ and $\xi_n^{hk}\in\partial J(Mu_n^{kh})$.

Note that the hypothesis $H(J)$(ii) and
$u\in H^1(0,T;V)$ imply that the sequence $\{\|\xi_{n}\|_{{X}^*}\}$ is uniformly bounded.
It follows from Lemma~\ref{lem3.2} that $\{\|\xi_n^{hk}\|_{{X}^*}\}$ is uniformly
bounded as well. Hence, we have
\begin{equation}\label{equation33}
\langle\xi_n-\xi_n^{hk},M(u_n-v^h_n)\rangle_{X^*\times X}\le C\|M(u_n-v^h_n)\|_X.
\end{equation}
Applying (\ref{equation31}) again, we obtain
\begin{equation}\label{equation34}
\langle \mathcal{R}_nu-\mathcal{R}_n^ku_n^{hk},
u_n-v^h_n\rangle\le C_2 k \bigg( 1+\sum_{j=1}^n\|e_j\|\bigg) \|u_n-v^h_n\|.
\end{equation}
Combining (\ref{equation32})--(\ref{equation34}) and applying the Cauchy inequality with $\varepsilon>0$,
we have
\begin{eqnarray}\label{equation35}
&&\|e_n\|^2_A-\|e_{n-1}\|^2_A
+2k(m_B-m_J\|M\|^2)\|e_n\|^2\le 2\langle Ae_n-Ae_{n-1},u_n-v_n^h\rangle\nonumber\\[2mm]
&&\quad
+\, C\, k\|u_n-v_n^h\|^2+\varepsilon k\|e_n\|^2+Ck^2\sum_{j=1}^{n-1}\|e_j\|^2
+Ck^3+C_2k^2\|e_n\|^2\nonumber\\
&&\qquad + \, C \, k \|M(u_n-v^h_n)\|_X + 2 k |S_n(v_n^h)|
+C k \|A\delta_n\|^2_{V^*}.
\end{eqnarray}

\noindent
Now we take
$\varepsilon=m_B-m_J\|M\|^2$ and $k_0=\frac{m_B-m_J\|M\|^2}{C_2}$,
which implies
\begin{equation*}
2\big(m_B-m_J\|M\|^2\big)-\varepsilon-kC_2>0,
\end{equation*}
for all $k<k_0$.
Subsequently, from (\ref{equation35}), we have
\begin{eqnarray}\label{equation36}
&&\|e_n\|^2_A-\|e_{n-1}\|^2_A\le 2\langle Ae_n-Ae_{n-1},u_n-v_n^h\rangle+Ck\|u_n-v_n^h\|^2+Ck^3
\nonumber\\
&&\quad + C k^2
\sum_{j=1}^{n-1}\|e_j\|^2+Ck\|M(u_n-v^h_n)\|_X+2k|S_n(v_n^h)|+Ck\|A\delta_n\|^2_{V^*}.
\end{eqnarray}
Now, we replace $n$ by $l$ in the above inequality,
and then sum it from $1$ to $n$,
where $1\le n\le N$ to get
\begin{eqnarray*}
&&\|e_n\|_A^2\le \|e_0\|_A^2+2\langle Ae_n,u_n-v_n^h\rangle+2\sum_{l=1}^{n-1}\langle Ae_l,(u_l-v_l^h)-(u_{l+1}-v_{l+1}^h)\rangle\nonumber\\
&&\quad+Ck\sum_{l=1}^n\bigg(\|u_l-v_l^h\|^2+\|M(u_l-v^h_l)\|_X+|S_l(v_l^h)|+\|A\delta_l\|^2_{V^*}\bigg)\nonumber\\
&&\qquad-2\langle Ae_0,u_1-v_1^h\rangle+Ck^2+Ck\sum_{l=1}^{n-1}\|e_l\|^2.
\end{eqnarray*}
This together with the following estimates
\begin{eqnarray*}
&&2\langle Ae_n,u_n-v_n^h\rangle\le \frac{1}{2}\|e_n\|_A^2+C\|u_n-v_n^h\|^2,\\ [2mm]
&&-2\langle Ae_0,u_1-v_1^h\rangle\le \|e_0\|_A^2+C\|u_1-v_1^h\|^2
\end{eqnarray*}
and
\begin{eqnarray*}
&&2\sum_{l=1}^{n-1}\langle Ae_l,(u_l-v_l^h)-(u_{l+1}-v_{l+1}^h)\rangle\le  2k\|A\|\sum_{l=1}^{n-1}\|e_l\|\|\delta(u_{l+1}-v_{l+1}^h)\|\nonumber\\
&&\quad\le C k\bigg(\sum_{l=1}^{n-1}\|e_l\|^2+\sum_{l=2}^n\|\delta(u_l-v_l^h)\|^2\bigg)
\end{eqnarray*}
implies that
\begin{eqnarray*}
&&\frac{1}{2}\|e_n\|_A^2\le 2\|e_0\|_A^2+C\|u_1-v_1^h\|^2+C\|u_n-v_n^h\|^2
+Ck\sum_{l=1}^{n-1}\|e_l\|^2+Ck^2\nonumber\\
&&\ +Ck\sum_{l=1}^n\bigg(\|\delta(u_l-v_l^h)\|^2
+\|u_l-v_l^h\|^2 + \|M(u_l-v^h_l)\|_X+|S_l(v_l^h)|
+\|A\delta_l\|^2_{V^*}\bigg).\nonumber\\
&&\qquad
\end{eqnarray*}
It follows from the discrete Gronwall inequality,  $H(A)$, and Lemma~\ref{gronwall}, that
\begin{eqnarray*}\label{cea}
&&\max_{0\le n\le N}\|e_n\|^2 \le C\bigg[k\sum_{l=1}^N\bigg(\|\delta(u_l-v_l^h)\|^2 + \|M(u_l-v^h_l)\|_X+|S_l(v_l^h)|+\|\delta_l\|^2\bigg)\nonumber\\
&&\hspace{5cm}
+ \|e_0\|^2+k^2+\max_{0\le n\le N}
\|u_n-v_n^h\|^2\bigg]
\end{eqnarray*}
for all $v_n^h\in V^h$.

We now summarize the results of the section in the form of a theorem.

\begin{Theorem}\label{CEA}
Suppose that assumptions of Lemma~\ref{lem3.1}
are satisfied.
Let  $u^{hk} \in V^h$ and $u \in H^1(0, T; V)$
be the solutions of Problems~\ref{fullydiscrete} and \ref{p1}, respectively.
Then, we have the estimate
\begin{eqnarray}\label{cea111}
&&
\hspace{-0.5cm}
\max_{0\le n\le N}\|u_n - u_n^{hk}\|^2 \le C\bigg[k\sum_{l=1}^N\bigg(\|\delta(u_l-v_l^h)\|^2 + \|M(u_l-v^h_l)\|_X+|S_l(v_l^h)|+\|\delta_l\|^2\bigg)
\nonumber
\\
&&\hspace{5cm}
+ \|e_0\|^2+k^2+\max_{0\le n\le N}
\|u_n-v_n^h\|^2\bigg]
\end{eqnarray}
for all $v_n^h\in V^h$.
\end{Theorem}

The inequality (\ref{cea111}) is called the C\'ea type inequality of the fully discrete approximation problem, Problem~\ref{fullydiscrete}.

\section{A quasistatic viscoelastic contact problem}
\label{Dynamic}
In this section we study the quasistatic contact problem between a viscoelastic body and a foundation.
The volume forces and surface tractions are supposed to
change slowly in time and therefore the acceleration
in the system is negligible.
Neglecting the inertial terms in the equation of motion leads to the quasistatic approximation for the process.
We show that the variational formulation of the quasistatic contact problem is a time-dependent hemivariational inequality in Problem~\ref{p1}.
For the latter, we apply the abstract result stated in Theorem~\ref{them3.2} and prove a result on existence and uniqueness of weak solution. Further, we use the fully discrete approximation method discussed in
Section~\ref{numerical} to study the numerical analysis of this contact problem and establish the result concerning optimal error estimate for the fully discrete scheme.

\subsection{Mathematical model
and its variational formulation}
\label{subsection1}

The physical setting of the contact problem is as follows.
A deformable viscoelastic body occupies an open bounded subset $\Omega$ of $\real^d$,
$d = 2$, $3$ in applications.
The volume forces of density $\fb_0$ act in $\Omega$
and surface tractions of density $\fb_N$ are applied on $\Gamma_2$. They both can depend on time.
We are interested in the quasi-static process of the mechanical state of the body on the time interval $[0, T]$ with $0< T < +\infty$.
The boundary $\Gamma = \partial \Omega$ of $\Omega$
is assumed to be
Lipschitz continuous and it consists of three
measurable parts $\Gamma_1$, $\Gamma_2$ and $\Gamma_3$
which are mutually disjoint,
and $m(\Gamma_1) > 0$. The unit outward
normal vector $\bnu$ exists a.e. on $\Gamma$.
We suppose that the body is clamped on part
$\Gamma_1$, and the body may come in contact
with an obstacle over the potential contact surface $\Gamma_3$. We also put
$Q = \Omega \times (0, T)$,
$\Sigma = \Gamma \times (0, T)$,
$\Sigma_1 = \Gamma_1 \times (0, T)$,
$\Sigma_2 = \Gamma_2 \times (0, T)$ and
$\Sigma_3 = \Gamma_3 \times (0, T)$.
We often do not indicate explicitly the dependence
of functions on the spatial variable $\bx \in \Omega$.

Let $\mathbb{S}^{d}$ denote the space of $d \times d$
symmetric matrices. The canonical inner
products and norms on $\real^d$ and $\mathbb{S}^{d}$
are given by
$$
\displaystyle
\bu \cdot \bv = u_i \, v_i, \quad \quad
\| \bv \| = (\bv \cdot \bv)^{1/2} \ \ \
{\rm for \ all} \ \ \bu, \bv \in \real^d,
$$
$$
\displaystyle
\bsigma : \btau = \sigma_{ij} \, \tau_{ij}, \quad \quad
\| \btau \| = (\btau : \btau )^{1/2} \ \ \
{\rm for \ all} \ \ \bsigma, \btau \in \mathbb{S}^{d}.
$$
\noindent
In what follows we always adopt the summation
convention over repeated indices.

Moreover, for a vector $\bxi \in \real^d$, the normal and tangential
components of $\bxi$ on the boundary are denoted by
$\xi_\nu = \bxi \cdot \bnu$ and $\bxi_\tau = \bxi - \xi_\nu \bnu$, respectively.
The normal and tangential components of the matrix
$\bsigma \in \mathbb{S}^{d}$ are defined on boundary
by $\sigma_\nu = (\bsigma \bnu) \cdot \bnu$ and
$\bsigma_\tau = \bsigma \bnu - \sigma_\nu \bnu$, respectively.

We denote by $\bu \colon Q \to \real^d$ the displacement vector,
by $\bsigma \colon Q \to \mathbb{S}^{d}$ the stress tensor and
by $\bvarepsilon (\bu) = (\varepsilon_{ij}(\bu))$ the linearized (small)
strain tensor, where $i$, $j = 1,\ldots, d$.
Recall that the components of the linearized strain tensor
are given by $\bvarepsilon (\bu) = 1/2 (u_{i,j} + u_{j,i})$, where
$u_{i,j} = \partial u_i / \partial x_j$.

The classical formulation of the contact problem reads as follows.

\medskip

\noindent
{\bf Problem} ${\mathcal P}$. {\it Find a displacement field
$\bu \colon Q \to\mathbb{R}^d$ and a stress field
$\bsigma \colon Q\to\mathbb{S}^d$ such that}
\begin{align}
&\Div \bsigma(t) + \fb_0(t)=0 \quad &\mbox{in}&\ Q,
\label{MOS2.1}\\[2mm]
&{\bsigma}(t) = {\mathscr A}\bvarepsilon({\bu}'(t))
+
{\mathscr B} \bvarepsilon(\bu(t))
+
\int_0^t{\mathscr C}(t-s) \bvarepsilon(\bu(s)) \,ds
\quad &\mbox{in}&\ Q,\label{MOS2.2}\\[2mm]
&\bu(t) = {\bzero}\quad &\mbox{on}&\ \Sigma_1,
\label{MOS2.3}\\[2mm]
&\bsigma(t)\bnu = \fb_N(t) \quad &\mbox{on}&\ \Sigma_2,
\label{MOS2.4}\\[2mm]
&-\sigma_\nu(t) \in \partial j_\nu(u_\nu(t)) \quad
&\mbox{on}&\ \Sigma_3,
\label{MOS2.5}\\[2mm]
&-\bsigma_\tau(t)\in \partial j_\tau ({\bu}_\tau(t)) \quad
&\mbox{on}&\ \Sigma_3, \label{MOS2.6}\\[2mm]
&\bu(0)=\bu_0 \quad &\mbox{in}&\ \Omega.
\label{MOS2.7}
\end{align}

\medskip

\noindent
The relation (\ref{MOS2.1}) represents the equilibrium equation
in which ``Div" denotes the divergence operator for tensor valued functions defined by
${\rm Div} \bsigma = (\sigma_{ij,j})$.
Equation (\ref{MOS2.2}) is the viscoelastic constitutive law with long memory,
where ${\mathscr A}$ and ${\mathscr B}$ are linear
viscosity and elasticity operators,
and ${\mathscr C}$ denotes the relaxation
operator.
Next, conditions (\ref{MOS2.3}) and
(\ref{MOS2.4}) represent the displacement and the traction boundary conditions.
The multivalued relations (\ref{MOS2.5}) and (\ref{MOS2.6})
are the contact and friction conditions, respectively, in
which
$\partial j_\nu$ and $\partial j_\tau$ denote the Clarke generalized gradients of prescribed locally Lipschitz
functions $j_\nu$ and $j_\tau$.
Finally, condition (\ref{MOS2.7}) represents the initial condition where $\bu_0$ denotes the initial displacement.
For concrete examples of boundary conditions (\ref{MOS2.5}) and (\ref{MOS2.6}), we refer to \cite{DL,HS,smo1,NP,PANA1,PANA7}.

Subsequently we introduce the spaces needed for the variational formulation.
Let $V$ be a closed subspace of $H^1(\Omega; \real^d)$
defined by
\begin{equation}\label{spaceV}
V = \{ \, \bv \in H^1(\Omega; \real^d) \mid \bv = 0 \ \ {\rm on} \ \Gamma_1 \, \}
\end{equation}

\noindent
and $H = L^2(\Omega; \real^d)$.
Then $(V, H, V^*)$ forms an evolution triple of spaces.
Moreover, the trace operator
is denoted by
$\gamma \colon V \to L^2(\Gamma; \real^d)$.
Given an element $\bv \in V$ we use the same notation
$\bv$ for the trace of $\bv$ on the boundary.
The space $V$ is equipped with the inner product and the corresponding norm given by
$$
\displaystyle
\langle \bu, \bv \rangle_V =
\langle \bvarepsilon (\bu), \bvarepsilon (\bv) \rangle_{\mathcal H},
\quad
\| \bv \| = \| \bvarepsilon (\bv) \|_{\mathcal H} \ \ {\rm for} \ \ \bu, \bv \in V,
$$

\noindent
where ${\mathcal H} = L^2(\Omega; \mathbb{S}^{d})$.
Since $m(\Gamma_1) > 0$,
from the Korn inequality
$\| \bv \|_{H^1(\Omega;\real^d)} \le c
\| \bvarepsilon (\bv) \|_{\mathcal H}$
for $\bv \in V$ with $c > 0$, it follows that $\| \cdot \|_{H^1(\Omega;\real^d)}$
and $\| \cdot \|$ are equivalent norms on $V$. In addition, we denote by $\mathcal{Q}_{\infty}$ the space of fourth-order tensor fields given by
\begin{equation*}
\mathcal{Q}_{\infty}
=\{\, \mathcal{E}=(\mathcal{E}_{ijkl}) \mid \mathcal{E}_{ijkl}=\mathcal{E}_{jikl}
=\mathcal{E}_{klij}\in L^\infty(\Omega),~1\le i,j,k,l\le d \, \}.
\end{equation*}

\medskip

We assume that the viscosity and elasticity tensors have the usual properties of ellipticity and symmetry.

\medskip

\lista{
	\item[$\underline{H({\mathscr A})}:$ ] \
	${\mathscr A} \colon \Omega \times \mathbb{S}^{d} \to \mathbb{S}^{d}$ is a viscosity tensor,
	${\mathscr A} = (a_{ijkl})\in\mathcal{Q}_{\infty}$ such that
	there exists $m_1 > 0$ satisfying
	${\mathscr A} \tau \cdot \tau \ge m_1
	\| \tau \|^2_{\mathbb{S}^{d}}$
	for all $\tau \in \mathbb{S}^{d}$, a.e. in $\Omega$.
}

\smallskip
\smallskip

\lista{
	\item[$\underline{H({\mathscr B})}:$] \
	${\mathscr B} \colon \Omega \times \mathbb{S}^{d} \to \mathbb{S}^{d}$ is an elasticity tensor,
	${\mathscr B} = (b_{ijkl})\in\mathcal{Q}_{\infty}$ such that there exists $m_2 > 0$ satisfying
	${\mathscr B} \tau \cdot \tau \ge m_2
	\| \tau \|^2_{\mathbb{S}^{d}}$
	for all $\tau \in \mathbb{S}^{d}$, a.e. in $\Omega$.
}

\smallskip
\smallskip

\lista{
	\item[$\underline{H({\mathscr C})}:$] \
	${\mathscr C} \colon [0,T] \to \mathcal{Q}_{\infty}$
	is Lipschitz continuous with constant $L_\mathscr{C}>0$.
}

\smallskip
\smallskip

The body forces, surface tractions and initial displacement satisfy

\smallskip
\smallskip

\lista{
	\item[$\underline{H(f)}:$]
	$\fb_0 \in L^2(0, T; L^2(\Omega; \real^d))$,
	$\fb_N \in L^2(0, T; L^2(\Gamma_2; \real^d))$,
	$\bu_0 \in V$.
}

\smallskip
\smallskip

The superpotentials satisfy

\smallskip
\smallskip

\noindent
$\underline{H(j_\nu)}:$ \quad
$j_\nu \colon \Gamma_3 \times \real \to \real$
is a function such that

\lista{
	\item[(i)]
	$j_\nu(\cdot, r)$ is measurable for all $r \in \real$,
	$j_\nu(\cdot, 0) \in L^1(\Gamma_3)$,
	\vspace{1mm}
	\item[(ii)]
	$j_\nu(\bx,\cdot)$ is locally Lipschitz for a.e.
	$\bx \in \Gamma_3$,
	\vspace{1mm}
	\item[(iii)]
	$| \partial j_\nu(\bx, r) | \le c_\nu (1+ | r |)$
	for a.e. $\bx \in \Gamma_3$, all $r \in \real$ with $c_\nu > 0$,
	\vspace{1mm}
    \item[(iv)] $(\eta_1-\eta_2)(r_1-r_2)\ge -m_\nu|r_1-r_2|^2$ for all $\eta_i\in \partial j_\nu(\bx,r_i)$, $r_i\in\real,~i=1,2$ for a.e. $\bx\in\Gamma_3$ with $m_\nu>0$.
}

\smallskip
\smallskip

\noindent
$\underline{H(j_\tau)}:$ \quad
$j_\tau \colon \Gamma_3 \times \real^d \to \real$
is a function such that

\lista{
	\item[(i)]
	$j_\tau(\cdot, \bxi)$ is measurable for all
	$\bxi \in \real^d$,
	$j_\tau(\cdot, {\bf0}) \in L^1(\Gamma_3)$,
	\vspace{1mm}
	\item[(ii)]
	$j_\tau(\bx, \cdot)$ is locally Lipschitz for a.e.
	$\bx \in \Gamma_3$,
	\vspace{1mm}
	\item[(iii)]
	$\| \partial j_\tau(\bx,\bxi) \|_{\real^d}
	\le c_\tau (1+ \| \bxi \|_{\real^d})$
	for a.e. $\bx \in \Gamma_3$, all $\bxi \in \real^d$ with $c_\tau > 0$,
	\vspace{1mm}
    \item[(iv)] $(\etab_1-\etab_2)\cdot(\bxi_1-\bxi_2)\ge -m_\tau\|\bxi_1-\bxi_2\|^2$ for all $\etab_i\in \partial j_\tau(\bx,\bxi_i)$, $\bxi_i\in \real^d,~i=1,2$ for a.e. $\bx\in\Gamma_3$ with $m_\tau>0$.
}

\smallskip
\smallskip

\noindent
In the hypotheses $H(j_\nu)$ and $H(j_\tau)$ the subdifferential is taken with respect
to the last variables of $j_\nu$ and $j_\tau$, respectively.

\medskip

Next, we define the operators
$A$, $B \in {\mathcal L}(V, V^*)$
by
\begin{equation}\label{333}
\langle A \bu, \bv \rangle_{V^* \times V} =
\langle {\mathscr A} \bvarepsilon(\bu),
\bvarepsilon (\bv) \rangle_{\mathcal H},
\quad
\langle B \bu, \bv \rangle_{V^* \times V} =
\langle {\mathscr B} \bvarepsilon(\bu),
\bvarepsilon (\bv) \rangle_{\mathcal H}
\end{equation}

\noindent
for $\bu$, $\bv \in V$, and the operator
${\mathcal R}\colon {\mathcal V} \to {\mathcal V^*}$
by
\begin{eqnarray}\label{definitionR}
\langle ({\mathcal R} \bw)(t), \bv \rangle_{V^* \times V}
= \Big\langle \int_0^t {\mathscr C}(t - s) \bvarepsilon (\bw(s)) \, ds, \bvarepsilon (\bv) \Big\rangle_{\mathcal H}
\end{eqnarray}
for all $\bw \in {\mathcal V}$, $\bv \in V$,
a.e. $t \in (0, T)$.

\medskip

To obtain the weak formulation of the problem (\ref{MOS2.1})--(\ref{MOS2.7}), we assume
the sufficient smoothness of the functions
involved, use the equilibrium
equation (\ref{MOS2.1}) and the Green formula.
We obtain
$$
\displaystyle
\langle \bsigma (t), \bvarepsilon (\bv))_{\mathcal H}
=
\langle \fb_0(t), \bv \rangle_H
+ \int_{\Gamma} \bsigma (t) \bnu \cdot \bv \,d\Gamma
$$

\noindent
for $\bv \in V$.
Taking into account the boundary condition (\ref{MOS2.3})
and (\ref{MOS2.4}),
we have
\begin{equation}\label{NEW0}
\langle \bsigma (t), \bvarepsilon (\bv)\rangle_{\mathcal H}
- \int_{\Gamma_3} \bsigma (t) \bnu \cdot \bv \,d\Gamma =
\langle \fb(t), \bv \rangle ,
\end{equation}

\noindent
where $\fb \in {\mathcal V}^*$ is given by
$\langle \fb(t), \bv \rangle =
\langle \fb_0(t), \bv \rangle_H
+ \langle \fb_N(t), \bv \rangle_{L^2(\Gamma_2; \real^d)}$
for $\bv \in V$.
On the other hand, by the ortogonality relation, cf. (6.33) in~\cite{smo1}, we get
\begin{equation}\label{NEW1}
\int_{\Gamma_3} \bsigma (t) \bnu \cdot \bv \,d\Gamma
=
\int_{\Gamma_3} \left( \sigma_\nu (t) v_\nu + \bsigma_\tau(t) \cdot \bv_\tau \right) \,d\Gamma.
\end{equation}

\noindent
The contact and friction boundary conditions
(\ref{MOS2.5}) and (\ref{MOS2.6}) can be equivalently formulated as follows
\begin{equation}\label{NEW2}
- \sigma_\nu(t) r \le j_\nu^0 (u_\nu; r) \ \ \mbox{for all} \ \ r \in \real, \ \ \
- \bsigma_\tau(t) \cdot \bxi \le j_\tau^0 (\bu_\tau; \bxi) \ \ \mbox{for all} \ \ \bxi \in \real^d.
\end{equation}

\noindent
Using (\ref{MOS2.2}), (\ref{333}), (\ref{NEW1}) and (\ref{NEW2}), from (\ref{NEW0}),
we obtain the following hemivariational inequality
which is a weak formulation of the problem
(\ref{MOS2.1})--(\ref{MOS2.7}):
find $\bu \colon (0, T) \to V$ such that
$\bu$, $\bu' \in {\mathcal V}$
and
\begin{eqnarray}\label{HVI}
\left\{
\begin{array}{llll}
\displaystyle
\langle A \bu'(t) + B \bu(t) + ({\mathcal R} \bu)(t), \bv \rangle +
\int_{\Gamma_3} \left( j_\nu^0 (u_\nu; v_\nu)
+ j_\tau^0 (\bu_\tau; \bv_\tau) \right)
\, d\Gamma \\ [3mm]
\qquad\qquad\qquad
\ge \langle \fb(t), \bv \rangle \ \
\mbox{for all} \ \ \bv \in V, \ \
\mbox{a.e.} \ \ t \in (0, T),
\\ [2mm]
\bu(0) = \bu_0.
\end{array}\right.
\end{eqnarray}

\subsection{Existence and uniqueness for contact
	problem}
\noindent
Let $X = L^2(\Gamma_3;\real^d)$ and
consider the functional
$J \colon X \to \real$ defined by
\begin{equation}\label{NEW3}
J(v) = \int_{\Gamma_3} \left(
j_\nu (\bx, v_\nu (\bx)) + j_\tau (\bx, \bv_\tau (\bx)) \right) \, d\Gamma \ \ \mbox{for all} \ \ \bv \in X.
\end{equation}

\noindent
Following~\cite[Theorem~5.1]{MOS1} and \cite[Corollary 4.15]{smo1},
we recall the following properties of the functional $J$.
\begin{Lemma}\label{LemmaJ}
	Under the hypotheses $H(j_\nu)$ and $H(j_\tau)$,  if, in
    addition,
    \begin{eqnarray}\label{clarkeregular}
    \left\{\begin{array}{lll}
    \mbox{either }j_\nu(\bx,\cdot)\mbox{ or }-j_\nu(\bx,\cdot)\mbox{ is regular and}\\[2mm]
    \mbox{either }j_\tau(\bx,\cdot)\mbox{ or }-j_\tau(\bx,\cdot)\mbox{ is regular,}
    \end{array}\right.
    \end{eqnarray}
    then the functional $J$
	defined by (\ref{NEW3}) satisfies
	
	\smallskip
	\lista{
		\item[\rm (i)]
		$J$ is Lipschitz continuous on bounded subsets of $X$,
		\smallskip
		\item[\rm (ii)]
		$\| \partial J(\bv) \|_{X^*}
		\le c_1 \left( 1 + \| \bv \|_X \right)$
		for all $\bv \in X$ with $c_1 = \max\{c_\tau,c_\nu\}$,
		\smallskip
        \item[\rm (iii)]
		for all $\bv$, $\bw \in X$, $\bxi\in\partial J(\bv)$ and $\etab\in\partial J(\bw)$, we have
	\begin{equation}\label{NEW5}
\langle \bxi-\etab,\bv-\bw\rangle_{X^*\times X}\ge -m_3\|\bv-\bw\|_X^2
	\end{equation}
with $m_3=m_\nu+m_\tau$,
		\smallskip
		\item[\rm (iv)]
		for all $\bv$, $\bw \in X$, we have
	}
	\begin{equation}\label{NEW4}
	J^0(\bv; \bw) = \int_{\Gamma_3} \left(
	j^0_\nu (v_\nu; w_\nu) + j^0_\tau (\bv_\tau; \bw_\tau) \right) \, d\Gamma
	\end{equation}

	\smallskip
	\lista{
		\item[]
		where $J^0(\bv; \bw)$ denotes the directional derivative of $J$
		at a point $\bv \in X$ in the direction
		$\bw \in X$.
	}
\end{Lemma}

Under our notation we associate with the hemivariational inequality (\ref{HVI}),
the following inclusion:
find $\bu \in {\mathcal V}$ such that
$\bu' \in {\mathcal V}$ and
\begin{eqnarray}\label{NEWINC11}
\left\{
\begin{array}{llll}
\displaystyle
\langle A \bu'(t) + B \bu(t) + ({\mathcal R} \bu)(t)
- \fb(t), \bv \rangle + J^0 (\gamma \bu(t); \gamma \bv) \ge 0 \\ [2mm]
\hspace{5.0cm} \mbox{for all } \  \bv \in V,
\ \mbox{a.e.}\ t \in (0, T), \\ [2mm]
\bu(0) = \bu_0.
\end{array}\right.
\end{eqnarray}

Note that if the hypotheses $H(j_\nu)$ and $H(j_\tau)$
hold,
then every solution to (\ref{NEWINC11}) is a solution to (\ref{HVI}).
The converse holds provided $j_\nu$ and $j_\tau$ satisfy the regularity condition (\ref{clarkeregular}).
These facts follow from the definition of the Clarke generalized gradient and Lemma~\ref{LemmaJ}.

The existence, uniqueness and regularity result for the hemivariational inequality (\ref{HVI}) is given in the following result.
\begin{Theorem}\label{existencetheorem}
	If the hypotheses $H({\mathscr A})$, $H({\mathscr B})$,
	 $H({\mathscr C})$,
	 $H(f)$, $H(j_\nu)$, $H(j_\tau)$, regularity
condition (\ref{clarkeregular}) hold, and the inequality $m_2> (m_\nu + m_\tau) \|\gamma\|^2$ is satisfied,
then problem (\ref{HVI})
has a unique solution $\bu \in H^1(0,T;V)$.
\end{Theorem}

\noindent
{\bf Proof.}
It follows from $H({\mathscr A})$ and $H({\mathscr B})$ that the operators $A$ and $B$
defined by (\ref{333}) satisfy $H(A)$ with $m_A = m_1$ and $H(B)$ with $m_B = m_2$, respectively. It is obvious from the definition of $\mathcal{R}$ (see (\ref{definitionR})) and hypothesis $H({\mathscr C})$ that $H(E)$ and $H(q)$ are satisfied with $E=I$ and $q=\mathscr{C}$.
Moreover, we put
$M = \gamma \in {\mathcal L}(V, X)$,
$\gamma$ is the trace operator.
It is a consequence of Lemma \ref{LemmaJ} that the functional $J$
given by (\ref{NEW3}) satisfies $H(J)$ with $c_J=c_1$ and $m_J=m_3$ (see Lemma \ref{LemmaJ}). Also $H(M)$ follows easily by the properties
of the trace operator. The conclusion is a consequence of Theorem~\ref{them3.2}, which completes the proof of this theorem.
\hfill$\Box$

\medskip

We say that a couple of functions $(\bu, \bsigma)$ which satisfies (\ref{MOS2.2}) and (\ref{HVI}) is called a weak solution to Problem~$\mathcal{P}$.
We
conclude that, under the assumptions of Theorem~\ref{existencetheorem}, Problem~$\mathcal{P}$ has a unique weak solution. Moreover, the weak
solution has the following regularity
$\bu \in H^1(0,T;V)$, $\bsigma \in L^2(0, T; L^2(\Omega,\mathbb{S}^{d}))$ and
${\rm Div} \, \bsigma \in {\mathcal V^*}$.

\subsection{Numerical analysis of contact problem}

In this section, we will apply the results from
Section~\ref{numerical} to establish an optimal order error estimate for the fully discrete solution of the contact problem in Problem~$\mathcal{P}$.
Here, we consider the frictionless boundary condition
on $\Gamma_3$, i.e., the frictional boundary (\ref{MOS2.6}) will be reduced to
\begin{equation}\label{frictionless}
\bsigma_\tau(t)={\bf 0} \ \ \mbox{on}\ \ \Sigma_3.
\end{equation}
In addition, without loss of generality, we may assume that $\bu_0={\bf0}$. We use the same the spaces as introduced in Section~\ref{subsection1}.
Then, consider the trace operator
$\gamma \colon V\to L^2(\Gamma_3;\real^d)$.
It follows from the Sobolev trace theorem that
\begin{eqnarray}\label{sobolevtrace}
\|\gamma\bv\|_{L^2(\Gamma_3;\real^d)}\le c_0\|\bv\|_{V}\ \,\mbox{ for all }\ \bv\in V
\end{eqnarray}
for some constant $c_0>0$, which depends only on  $\Omega$, $\Gamma_1$ and $\Gamma_3$.
Let $X=L^2(\Gamma_3)$ and define the operators $\gamma_\nu\colon L^2(\Gamma_3;\real^d)\to X$, $\gamma_\nu\bv = v_\nu$ for
$\bv\in L^2(\Gamma_3;\real^d)$,
and $M = \gamma_\nu\circ \gamma\colon V\to X$.
We also consider the functional
$J \colon X \to \real$ defined by
\begin{equation*}
J(v) = \int_{\Gamma_3}
j_\nu (\bx, v_\nu (\bx))  \,d\,\Gamma \ \ \mbox{for all} \ \ v \in X.
\end{equation*}
If either $j_\nu(\bx,\cdot)$  or $-j_\nu(\bx,\cdot)$ is regular and $H(j_\nu)$ holds, then by Lemma \ref{LemmaJ}, (\ref{333}) and (\ref{definitionR}), the contact problem (\ref{MOS2.1})--(\ref{MOS2.7}) with $j_\tau=0$ has following equivalent variational formulation.
\begin{Problem}\label{variationalformulation}
Find $\bu \in \mathcal{V}$ such that $\bu'\in\mathcal{V}$  and
\begin{eqnarray}\label{NEWINC12}
\left\{
\begin{array}{llll}
\displaystyle
\langle A \bu'(t) + B \bu(t) + ({\mathcal R} \bu)(t)
- \fb(t), \bv-\bu(t) \rangle \\ [2mm]
\qquad
+ \, J^0 (\gamma \bu(t); \gamma \bv-\gamma \bu(t))
\ge 0 \ \ \mbox{for all} \ \ \bv \in V,
\ \mbox{a.e.}\ t \in (0, T),  \\ [2mm]
\bu(0) = {\bf0}.
\end{array}\right.
\end{eqnarray}
\end{Problem}

From Theorem~\ref{existencetheorem}, we deduce that under the hypotheses  $H({\mathscr A})$, $H({\mathscr B})$,
	 $H({\mathscr C})$, and
	 $H(f)$, $H(j_\nu)$. If either $j_\nu(\bx,\cdot)$  or $-j_\nu(\bx,\cdot)$ is regular and the inequality
$m_2> (m_\nu + m_\tau) \|\gamma\|^2$ hold, then Problem~\ref{variationalformulation} has a unique solution $\bu \in H^1(0, T; V)$.

Next, we pass to the numerical approximation of
Problem~\ref{variationalformulation}.
Likewise in Section~\ref{numerical}, for an integer $N>0$, let $k=\frac{T}{N}$ be the time step length.
For simplicity, we suppose that $\Omega$ is a polygonal/polyhedral domain and express the three parts
of the boundary, $\Gamma_k$, $k=1$, $2$, $3$,
as a union of closed flat components
with disjoint interiors
\begin{equation*}
\overline{\Gamma}_j
=\cup_{i=1}^{i_j}\Gamma_{j,i},\ \ \,1\le j\le 3.
\end{equation*}
Subsequently, we consider a regular family of meshes $\{\mathcal T^h\}$ that partition ${\bar{\Omega}}$ into triangles/tetrahedrons compatible
with the splitting of the boundary $\partial \Omega$
into $\Gamma_{j,i}$, $1\le i\le i_j$, $1\le j\le 3$.
This means that if the intersection of one side/face of an element with one set $\Gamma_{j,i}$ has a positive
measure with respect to $\Gamma_{j,i}$,
then the side/face lies entirely in $\Gamma_{j,i}$.
Corresponding to the family $\{\mathcal T^h\}$,
we define the linear element space
\begin{equation*}
V^h
=\big\{\, \bv^h\in C(\overline{\Omega}; \real^d) \mid
\bv^h|_U\in \mathbb{P}_1(U)^d,~U\in\mathcal{T}^h,~\bv^h=0\mbox{ on }\Gamma_1\, \big\},
\end{equation*}
where $\mathbb{P}_1(U)^d$ denotes a set of all linear functions whose domain of definition is $U$ (cf.~\cite[p. 70]{HMP}).

Now, we are in a position to formulate the following fully discrete approximation problem for Problem~\ref{variationalformulation}.
\begin{Problem}\label{fullydiscreteproblem}
Find
$\bu^{hk}=\{\bu_n^{hk}\}\subset V^h$
such that $\bu_0^{hk}={\bf0}$ and
\begin{eqnarray}\label{fullydiscretesequence}
\left\{
\begin{array}{llll}
\displaystyle
\langle A\delta \bu_n^{hk}+B \bu_n^{hk}+\mathcal{R}_n^{k}\bu^{hk},\bv^h-\bu_n^{hk}\rangle +J^0(M\bu_n^{hk};M\bv^{h}-M\bu_n^{hk})
\\[2mm]
\hspace{3cm}\ge\langle \fb_n,\bv^h-\bu_n^{hk}\rangle
\ \ \mbox{for all} \ \ \bv^h\in V^h
\end{array}\right.
\end{eqnarray}
for all $n=1,2,\ldots,N$.
\end{Problem}

In the sequel, we assume that the solution of
Problem~\ref{variationalformulation} has the following additional regularity
\begin{eqnarray}\label{regularsolution}
\left\{\begin{array}{lll}
\bu\in H^1(0,T;H^2(\Omega)), \ \, \bu''\in L^2(0,T;V),\\[2mm]
u_\nu|_{\Gamma_{3,i}}\in C(0,T;H^2(\Gamma_{3,i})),\ \, \sigma_\nu|_{\Gamma_{3,i}}\in C(0,T;L^2(\Gamma_{3,i}))
\end{array}\right.
\end{eqnarray}
for $1\le i\le i_3$. Then the function
$(t,\bx)\to \bu(t,\bx)$ is continuous.
This means that the pointwise values of $\bu$ are well-defined.
So, take $\bv_n^h=\Pi^h\bu_n\in V^h$ to be the finite element interpolant of $\bu_n(\bx) =\bu(t_n,\bx)$,
where
$\Pi^h\bu_n$ denotes the piecewise constant Lagrange interpolation of $\bu_n$ (cf.~\cite[p. 122]{HMP}).
We use the C\'ea type inequality (\ref{cea}) to get
\begin{eqnarray}
&&
\hspace{-0.5cm}
\max_{1\le n\le N}\|\bu_n-\bu_n^{hk}\|^2_V \le C\bigg[k^2+\max_{1\le n\le N}\|\bu_n-\Pi^h\bu_n\|^2_V+k\sum_{n=1}^N\bigg(\|\delta(\bu_n-\Pi^h\bu_n)\|^2_V\nonumber\\
&&\quad + \|u_{n,\nu}-\Pi^hu_{n,\nu}\|_X+|S_l(\Pi^h\bu_n)|+\|\bdelta_n\|^2_V\bigg)\bigg],
\end{eqnarray}
where
\begin{eqnarray*}
&&\bdelta_n =\delta\bu_n-\bu_n',\\[2mm]
&&S_n(\bv) =\langle A\bu_n'+B\bu_n+\mathcal{R}_n\bu-\fb_n,\bv-\bu_n\rangle+J^0(M\bu_n;M\bv-M\bu_n).
\end{eqnarray*}

It follows from \cite[Lemma 11.5]{HR} that
\begin{equation*}
\|\bdelta_n\|_V\le \|\bu''\|_{L^1(t_{n-1},t_n;V)}.
\end{equation*}
This together with H$\ddot{\rm{o}}$lder inequality implies that
\begin{equation*}
\|\bdelta_n\|_V^2\le k\|\bu''\|_{L^2(t_{n-1},t_n;V)}^2
\end{equation*}
and
\begin{equation}\label{bdeltaestimate}
k\sum_{n=1}^N\|\bdelta_n\|^2_V\le k^2\|\bu''\|_{L^2(0,T;V)}^2.
\end{equation}
Next, we use the fact
\begin{equation*}
\delta(\bu_n-\Pi^h\bu_n)
=\frac{1}{k}\int_{t_{n-1}}^{t_n}\big(\bu'(s)-\Pi^h\bu'(s)\big)\, ds
\end{equation*}
to obtain
\begin{equation*}
\|\delta(\bu_n-\Pi^h\bu_n)\|_V^2\le
\frac{1}{k}\int_{t_{n-1}}^{t_n}\big\|\bu'(s)-\Pi^h\bu'(s)\big\|^2_V\, ds.
\end{equation*}
Hence, we have
\begin{equation*}
k\sum_{n=1}^N\|\delta(\bu_n-\Pi^h\bu_n)\|_V^2\le \int_0^T\big\|\bu'(s)-\Pi^h\bu'(s)\big\|^2_V\, ds
\end{equation*}
and
\begin{equation}\label{deltaestimate}
k\sum_{n=1}^N\|\delta(\bu_n-\Pi^h\bu_n)\|_V^2\le Ch^2\|\bu'\|_{L^2(0,T;V)}^2.
\end{equation}
Recall that $\Pi^hu_{n,\nu}$ is the finite element interpolant of $u_{n,\nu}$ on each component $\Gamma_{3,i}$.
Combining (\ref{deltaestimate}) with the hypothesis (\ref{regularsolution}), we get
\begin{equation}\label{Xestimate}
k\sum_{n=1}^N\|u_{n,\nu}-\Pi^hu_{n,\nu}\|_X
\le Ch^2\sum_{i=1}^{i_3}
\|u_\nu\|_{L^\infty(0,T;H^2(\Gamma_{3,i}))}.
\end{equation}

On the other hand, we estimate the residual quantity $|S_n(\bv)|$. To this end, we use the fact (see (\ref{NEW0}), (\ref{NEW1}) and (\ref{frictionless})) that
\begin{equation*}
\langle A\bu_n'+B\bu_n+\mathcal{R}_n\bu-\fb_n,\bv\rangle
=\int_{\Gamma_3}\sigma_\nu(t)v_\nu \, d\Gamma
\ \ \mbox{for all}\ \ \bv\in V
\end{equation*}
to get
\begin{equation*}
S_n(\bv)
=\int_{\Gamma_3}\big(\sigma_{n,\nu}+\xi_n\big)
\big(\Pi^hu_{n,\nu}-u_{n,\nu}\big) \, d\Gamma
\end{equation*}
for some $\xi_n\in\partial J(M\bu_{n})$.
This implies
\begin{equation*}
|S_n(\Pi^h\bu_n)|\le C\|\Pi^hu_{n,\nu}-u_{n,\nu}\|_X,
\end{equation*}
and, therefore, we have
\begin{equation*}
k\sum_{n=1}^N|S_n(\Pi^h\bu_n)|\le Ch^2\sum_{i=1}^{i_3}
\|u_\nu\|_{L^\infty(0,T;H^2(\Gamma_{3,i}))}.
\end{equation*}
This estimate together with (\ref{bdeltaestimate})--(\ref{Xestimate}) implies
the following optimal estimate for the fully discrete scheme (\ref{fullydiscretesequence}).

\begin{Theorem}\label{Lasttheorem}
Assume that $\bu$ and $\bu^{hk}$ are solutions to Problems~\ref{variationalformulation} and~\ref{fullydiscreteproblem}, respectively,
and the regularity condition (\ref{regularsolution}) holds.
Then, we have
\begin{equation*}
\max_{1\le n\le N}\|\bu_n-\bu_n^{hk}\|_V
\le C(k+h),
\end{equation*}
where $C>0$ is independent of $k$ and $h$.
\end{Theorem}

In the optimal error estimate of Theorem~\ref{Lasttheorem},
the method is of first-order in spatial mesh size
and in the time step.

\end{document}